\documentclass[11pt]{amsart}
\usepackage[left=1.25in,right=1.25in,top=1in,bottom=1in]{geometry}
\usepackage{amsmath,amsfonts,amssymb,amsthm}
\usepackage[all]{xy}
\usepackage[shortlabels]{enumitem}
\usepackage{hyperref}
\usepackage{leftidx}
\usepackage{cleveref} 

\theoremstyle{plain}
\newtheorem{thm}{Theorem}[section]
\newtheorem{lemma}[thm]{Lemma}
\newtheorem{proposition}[thm]{Proposition}
\newtheorem{corollary}[thm]{Corollary}

\theoremstyle{definition}
\newtheorem{definition}[thm]{Definition}

\newtheorem{example}[thm]{Example}

\theoremstyle{remark}
\newtheorem{remark}[thm]{Remark}

\newcommand{\inj}{\hookrightarrow}

\newcommand{\Mfg}{\mathcal{M}_{\mathrm{fg}}}

\newcommand{\sma}{\wedge}

\newcommand{\Cts}{\mathrm{Maps}_{cts}}

\newcommand{\defcommand}[2]{
  \providecommand{#1}{}
  \renewcommand{#1}{#2}
}

\usepackage{etoolbox}
\newcommand{\define}[4]{\expandafter#1\csname#3#4\endcsname{#2{#4}}}
\forcsvlist{\define{\defcommand}{\mathbb}{}}{A,C,F,G,H,L,N,Q,R,T,Z}
\forcsvlist{\define{\defcommand}{\mathcal}{s}}{F,M,O}
\forcsvlist{\define{\defcommand}{\mathsf}{}}{Ab,Alg,Pro,Ind,Pipes,PipeRings,Sets,Rings,
  Mod,Cat,CLN,Comod,ComodAlg,Def,DefFrob,Top,Spaces,FGL,FG,Gpd,Groups,Ho,
  IsogMod,IsogAlg,StagedDef}
\forcsvlist{\define{\DeclareMathOperator}{}{}}{
  Ann,Aut,Der,End,Ext,Frob,Gal,Gr,Hom,
  height,Iso,Inf,Maps,rank,Spa,Spec,Spf,Tor,Tr}
\forcsvlist{\define{\defcommand}{}{}}{LE}

\newcommand{\CplMod}{\mathsf{Mod}^\wedge}

\newcommand{\Coh}{C}
\newcommand{\CLaur}{\Lambda}
\newcommand{\aug}{\mathrm{aug}}

\newcommand{\ol}[1]{\overline{#1}}

\newcommand{\Einf}{\mathcal{E}_\infty}

\DeclareMathOperator*{\holim}{holim}

\newcommand{\pullbackcorner}[1][dr]{\save*!/#1-2pc/#1:(-1,1)@^{|-}\restore}

\title{Localizations of Morava $E$-theory and deformations of formal groups}
\author{Paul VanKoughnett}

\begin{document}

\begin{abstract}
  We study the relationship between the transchromatic localizations of Morava
  $E$-theory, $L_{K(n-1)}E_n$, and formal groups. In particular, we show that
  the coefficient ring $\pi_0L_{K(n-1)}E_n$ has a modular interpretation,
  representing deformations of formal groups with certain extra structure, and
  derive similar descriptions of the cooperations algebra and $E_{n-1}$-homology
  of this spectrum. As an application, we show that $L_{K(1)}E_2$ has exotic
  $\Einf$ structures not obtained by $K(1)$-localizing the $\Einf$ ring $E_2$.
\end{abstract}

\maketitle

\setcounter{tocdepth}{1}
\tableofcontents

\section{Introduction}\label{introduction}

Stable homotopy theory is concerned with the study and construction of
cohomology theories. Many interesting cohomology theories are complex
orientable, meaning that they can be equipped with a system of Chern classes for
complex vector bundles; as it happens, this data is neatly encoded in a purely
algebraic object known as a formal group. The relation between formal groups and
stable homotopy theory in general has been an enduring theme over the last fifty
years. It was first sounded by Quillen, who proved in \cite{Quillen} that the
complex cobordism theory $MU$, which carries the universal theory of Chern
classes, also carries the universal formal group law. The Landweber exact
functor theorem \cite{Landweber}, in a deep baritone, showed that cohomology
theories could be constructed out of formal groups, and this version of the
theme has been trumpeted again by the modern theory of topological modular forms
\cite{TMF} and its higher generalizations \cite{BL}. Others like Morava and
Ravenel took up the call by relating the height filtration on formal groups to
the $K(n)$-local filtration on the stable homotopy category \cite{Morava,
  Rav}, and we now know that the classification of formal groups is linked to
deep structural properties of the stable homotopy category as a whole \cite{DHS,
  HSmith}.

The theme has found its clearest expression in the study of Morava $E$-theory.
This is a complex oriented cohomology theory constructed from the algebraic
geometry of deformations of formal groups. In \cite{LT}, Lubin and Tate proved
that a height $n$ formal group $\Gamma$ over a perfect, characteristic $p$ field
$k$ has a universal deformation which lives over the ring $Wk[[u_1, \dotsc,
u_{n-1}]]$. The parameters $u_i$ control the height of the deformation: for
example, inverting $u_{n-1}$ forces the height to be at most $n-1$. By the
theorem of Goerss, Hopkins, and Miller \cite{GH1}, there is an essentially
\textit{unique} complex oriented, $K(n)$-local, $\Einf$ ring spectrum, called
\textbf{Morava $E$-theory} $E = E(k,\Gamma)$, with
\[ \pi_*E = Wk[[u_1, \dotsc, u_{n-1}]][u^{\pm 1}], |u| = 2, |u_i| = 0, \] and
with formal group the universal deformation defined by Lubin and Tate.

This theorem suggests that the relationship between stable homotopy theory and
the algebraic geometry of formal groups is extremely close when localized at a
single prime and height. In particular, basically all topological facts about
$E$-theory should be expressible in terms of formal groups. For example:
\begin{enumerate}
\item The profinite group $\G_n$ of automorphisms of the field $k$ and the
  formal group $\Gamma$ acts on the Lubin-Tate ring, and this extends to an
  action of $\G_n$ on $E_n$ by $\Einf$ maps. By a theorem of Devinatz and
  Hopkins (\cite{DH}, and see \Cref{E cooperations} in this document),
  \[ \pi_*L_{K(n)}(E_n^{\sma (s+1)}) \cong \Hom_{cts}(\G_n^{\times s},
    \pi_*E_n). \] This means that the $K(n)$-local $E_n$-based Adams spectral
  sequence for the sphere takes the form
  \[ E_2^{st} = H^s(\G_n, \pi_tE_n) \Rightarrow \pi_{t-s}L_{K(n)}S. \] Moreover,
  $L_{K(n)}S$ is the homotopy fixed points of $\G_n$ acting on $E_n$, in a sense
  described by \cite{DH}. This also means that the $E_*$-comodule structure on
  the completed $E$-homology of a space or spectrum is just a continuous
  $\G_n$-action.
\item The completed $E$-homology of a $K(n)$-local $\Einf$ ring spectrum carries
  power operations, which are parametrized by isogenies of deformations of the
  formal group $\Gamma$ \cite{AHS, Rezk}.
\item The $E$-theory of Eilenberg-MacLane spaces for abelian groups can be
  described in terms of exterior powers of the universal deformation of $\Gamma$
  \cite{Peterson, HLamb}.
\item For nonabelian groups $G$, there is a character map from $E^0BG$ with
  image in a ring of `generalized class functions' on conjugacy classes of
  $n$-tuples of commuting elements of $p$-power order in $G$, which becomes an
  isomorphism after base changing to a ring parametrizing certain level
  structures on formal groups \cite{HKR, StapletonHKR}.
\end{enumerate}

Of the transchromatic objects straddling heights $n-1$ and $n$, one of the most
basic is the $K(n-1)$-localization of a height $n$ $E$-theory, known in this
document as $L_{K(n-1)}E_n$ or just $\LE$. This is even periodic, with
\[ \LE_0 = Wk[[u_1, \dotsc, u_{n-1}]][u_{n-1}^{-1}]^{\wedge}_{(p, u_1, \dotsc,
    u_{n-2})} \cong Wk((u_{n-1}))^\wedge_p[[u_1, \dotsc, u_{n-2}]], \] a
complete local ring with residue field $k((u_{n-1}))$, over which there is a
naturally defined height $n-1$ formal group $\H = \G^u \otimes k((u_{n-1}))$.
This looks very much like the height $n-1$ $E$-theory associated to $\H$ -- the
only problem being that the field $k((u_{n-1}))$ is not perfect, so that the
Lubin-Tate theorem does not apply.

This article has two goals. The first is to clarify the relationship between
$\LE$ and formal groups. As we prove, the ring $\LE_0$ classifies deformations
of $\H$ together with certain extra structure, which can be briefly described as
follows.

\begin{thm}[{\Cref{LE representability}}] Continuous ring homomorphisms $\LE_0
  \to R$, for a complete local ring $R$, correspond to deformations of $\H$ over
  $R$, together with a choice of last Lubin-Tate coordinate,
  \[ j:Wk((u_{n-1}))^\wedge_p \to R. \]
\end{thm}

We construct similar modular interpretations of other invariants of $\LE$,
namely the completed $E_{n-1}$-homology
\[ ( E_{n-1}^\wedge )_*\LE = \pi_*L_{K(n-1)}(E_{n-1} \sma \LE) =
  \pi_*L_{K(n-1)}(E_{n-1} \sma E_n) \] and the completed cooperations algebra
\[ \LE_*^\wedge \LE = \pi_*L_{K(n-1)}(E_n \sma \LE) = \pi_*L_{K(n-1)}(E_n \sma
  E_n). \] Like $\LE_*$, these rings are even periodic and it suffices to
describe them in degree 0.

\begin{thm}[{\Cref{LE cooperations}, \Cref{FE theorem}}] Continuous maps $(
  E_{n-1}^\wedge )_0\LE \to R$, for complete local rings $R$, represent pairs of
  a deformation of $\H$ and a deformation of the height $n-1$ formal group
  $\Gamma$ over $R$, together with an isomorphism of these formal groups over
  $R/\mathfrak{m}$, and a choice of last Lubin-Tate coordinate,
  \[ j:Wk((u_{n-1}))^\wedge_p \to R. \]

  Continuous maps $\LE_0^\wedge \LE \to R$, for complete local rings $R$,
  represent pairs of deformations of $\H$ over $R$, pairs of last Lubin-Tate
  coordinates
  \[ j, j':Wk((u_{n-1}))^\wedge_p \to R, \] and an isomorphism of the
  deformations over $R/\mathfrak{m}$.
\end{thm}

Interestingly, the ``last Lubin-Tate coordinate'' appearing in these theorems
exists mostly separately from the formal group deformation data, a fact which is
useful in doing constructions with these objects.

The second goal is to study power operations and $\Einf$ structures on the
spectra $\LE$. Here we restrict ourselves to $n = 2$ and $n-1 = 1$, so that we
can take advantage of the obstruction theory defined in \cite{GH2}, which
allows one to construct $K(1)$-local $\Einf$ ring spectra out of power
operations data known as a $\theta$-algebra structure. We prove:

\begin{thm}[{\Cref{LK(1)E2}, \Cref{LK(1)E2 nonisomorphic}}] There are $\Einf$
  structures on $L_{K(1)}E_2$ not obtained by $K(1)$-localizing the unique
  $\Einf$ structure on $E_2$.
\end{thm}

Besides serving as an instructive calculation in $K(1)$-local obstruction
theory, this result, and the methods used to prove it, point the way towards
studying the transchromatic behavior of highly structured ring spectra and their
power operations. In particular, these exotic $\Einf$-algebras are
$K(1)$-localizations of $K(2)$-local \textit{spectra}, and are
\textit{$\Einf$-algebras}, but only one of them -- the $K(1)$-localization of
the canonical $\Einf$-algebra $E_2$ -- is a $K(1)$-localization of a
$K(2)$-local $\Einf$-algebra. It is likely that, as is the case here,
$K(n)$-local power operations generally satisfy integrality conditions that can
be distorted on their $K(n-1)$-localizations. In future work, the author hopes
to prove the analogous statement at all heights $n$.

\subsection{Other approaches}\label{three continuities}
Before outlining the paper, let's highlight some of the difficulties afoot in
this transchromatic situation, and compare the approach taken here with some of
the other ones in the literature. Again, $E_0$ carries the universal deformation
of a height $n$ formal group $\Gamma$ over a finite field $k$, and roughly
speaking, passing to the $K(n-1)$-localization should be equivalent to
restricting to deformations which are height exactly $n-1$. But it's unclear how
exactly to think of these objects as deformations of $\Gamma$, since $\Gamma$
itself is not a deformation of this type. From an algebraic point of view, the
ring map $E_0 \to \LE_0$ is analogous to the map
\[ \Z_p[[x]] \to \Z_p((x))^\wedge_p. \] This map is not continuous with respect
to the maximal ideal topologies on the two rings, and so does not induce a map
\[ \Spf \Z_p[[x]] \leftarrow \Spf \Z_p((x))^\wedge_p; \] one is somehow trying
to take the complement of the unique point of $\Spf \Z_p[[x]]$.

One possible approach is to treat $\Z_p((x))^\wedge_p$ as a topological ring
such that the above map is continuous; that is, a basis of neighborhoods of 0 is
given by $(p, x)^r\Z_p[[x]]$. This is a topology with comparatively many open
sets. The subring $\Z_p((x))$ is open, and a continuous map $\Z_p((x)) \to R$ is
simply a continuous map $\Z_p[[x]] \to R$ sending $x$ to a unit; on the other
hand, $\Z_p((x))$ is not dense in $\Z_p((x))^\wedge_p$, so that some choices
with no simple description are needed to define maps $\Z_p((x))^\wedge_p \to R$.
Torii uses this topology to study localizations of $E$-theory \cite{Torii}.

Another approach is to treat the maximal ideal topology and the topology coming
from $\Z_p[[x]]$ as coexisting simultaneously on $\Z_p((x))^\wedge_p$. One way
to do this is to treat $\Z_p((x))^\wedge_p$ as a pro-object which is the limit
of the rings $\Z/p^n((x))$, each of which carries the $x$-adic topology. This
idea is explored by \cite{MGPS}, who develop a theory of ``pipe rings'' based on
ideas from the theory of higher local fields in number theory \cite{Morrow,
  Kato}, and prove analogous theorems about $\LE$ in this context. There are
some serious difficulties in the algebraic geometry that have so far presented
obstacles to more widespread adoption of this program.

We remark that some of the results here have pipe parallels: for example, just
as $\LE_0$ represents augmented deformations of the formal group $\H$ as a
complete local ring, it also does as a pipe ring. Additionally, we will use some
basic concepts from pipe theory to prove results about $\theta$-algebras in
section 6.

Third, we should mention the recently introduced framework of condensed
\cite{ClausenScholze} or pyknotic \cite{BarwickHaine} mathematics, which
replaces topological groups, rings, and so on with better-behaved categories of
sheaves on the site of compact Hausdorff spaces. The author believes that this
framework may present a way around some of the difficulties in the theory of
pipe rings, but this has not been worked out.

In this paper we take the fourth and simplest approach, which is to think of
$\LE_0$ as a topological ring with its maximal ideal topology (for example,
$\Z_p((x))^\wedge_p$ carries the $p$-adic topology). Thus a continuous map
$\Z_p((x))^\wedge_p \to R$, for some complete local ring $R$, descends to a map
$\F_p((x)) \to R/\mathfrak{m}$ which may not have any nice continuity
properties; note that $\F_p((x))$ has infinite transcendence degree over
$\F_p(x)$, so that such a map is not determined by the image of $x$. Of course,
the multitude of such maps is part of the reason for the multitude of $\Einf$
structures found in \Cref{LK(1)E2}.

Let us finally comment that one can \textit{also} view $\LE_0$ as carrying not a
formal group but a $p$-divisible group. This is the result of base changing the
$p$-power torsion subgroups of $\G^u$ to $\LE_0$; it fits into an exact sequence
\[ 0 \to \G^{for} \to \G \to \G^{et} \to 0 \] where $\G^{for}$ is the formal
group of $\LE$, and $\G^{et}$ is a Galois twist of a constant group scheme of
the form $\underline{\Q_p/\Z_p}$. Homotopy theorists have long thought that this
$p$-divisible group structure should manifest itself in transchromatic
phenomena, in particular due to the use of $p$-divisible groups of abelian
varieties in defining topological automorphic forms \cite{LurieSEC, LurieECI,
  LurieECII, BL}; the best work on this subject is still \cite{StapletonHKR},
which uses $p$-divisible groups like these to produce character maps.

Again, we take a simpler approach here, forgetting about the \'etale part of
this $p$-divisible group and only considering the deformation of the formal
group $\H$.

\subsection{Outline of the paper}
There are two background sections: section 2 concerns topics from algebraic
geometry (particularly concerning formal groups and deformation theory), and
section 3 discusses Morava $E$-theory. For purposes of comparison, we review the
Devinatz-Hopkins calculation (as \Cref{E cooperations}) of the completed
cooperations coalgebra $E_*^\wedge E$.

In section 4, we describe basic facts about $\LE$, and prove (\Cref{LE
  representability}) that its coefficient ring represents so-called augmented
deformations of a height $n-1$ formal group.

In section 5, we describe the cooperations $\LE_*^\wedge \LE$ (\Cref{LE
  cooperations}), and the height $n-1$ $E$-theory $(E_{n-1})^\wedge_*\LE$
(\Cref{FE theorem}).

The final section 6 concerns power operations and $\Einf$ structures, and works
specifically with $L_{K(1)}E_2$. After a detour into the theory of
$\theta$-algebras, the existence of exotic $\Einf$ structures on this spectrum
is proved in \Cref{LK(1)E2}.

\subsection{Acknowledgements}
This article is a chunk of my PhD thesis, and couldn't have been finished, let
alone started, without the help of my advisor, Paul Goerss. I can't thank Paul
enough for his wisdom, support, and insistence that I make my confusion specific
(a joke originally due to Steve Wilson).

I want to thank \"Ozg\"ur Bayindir, Irina Bobkova, Jack Morava, Eric Peterson,
Jay Shah, Joel Specter, Nat Stapleton, and Dylan Wilson for interesting
questions and helpful mathematical conversations. I want to especially thank
Irina for some much-needed pressure to put these results into readable form. I
want to thank Dominik Absmeier for pointing out errors in earlier versions of some of the proofs.

Finally, I want to thank those who kept me sane at one point or another while I
was writing the thesis: my friends, Pete Davis, Ilana Shalowitz, Sara Shaw, Joel
Specter (again), and Talia Lavin; my family, Diane, Hale, and Will VanKoughnett;
and my tree in the wind, Rox Sayde.

\subsection{Notation}
For the sake of reference, here is some notation introduced elsewhere in the
paper:

$k$ is a field of characteristic $p$, generally perfect, often even finite. $W$
denotes the Witt vectors functor. $\CLaur$ is the completed Laurent series ring
$Wk((u_{n-1}))^\wedge_p$.

$\CLN$ is the category of complete local noetherian rings and continuous ring
homomorphisms -- for $\CLN_A$, see \Cref{CLN}. $\Gpd$ is the category of
groupoids. Other categories are written in \textsf{sans-serif}, but typically
identifiable from their names.

$\Gamma_n$ (sometimes just $\Gamma$) is a formal group over $k$ of finite height $n$.

$E = E_n = E(k,\Gamma)$ is the Morava $E$-theory for $(k, \Gamma)$; this is even
periodic, with periodicity class $u$ in degree $2$ and $\pi_0 E$ non-canonically
isomorphic to $Wk[[u_1, \dotsc, u_{n-1}]]$. $\LE$ is the localization
$L_{K(n-1)}E_n$. At a certain point, $F$ is used for $E_{n-1}$, and $K$ for $E_1$.

$BP$ is the Brown-Peterson spectrum with $BP_* = \Z_{(p)}[v_1, v_2, \dotsc]$,
where $|v_i| = 2(p^i-1)$; the map $BP \to E$ sends $v_n \mapsto u^{p^n-1}$, $v_i
\mapsto 0$ for $i > n$, and and $v_i \mapsto u^{p^i-1}u_i$ for $i < n$. $I_n$ is
the ideal $(p, v_1, \dotsc, v_{n-1})$, or its image in any $BP_*$-module.

Notation like $E_*^\wedge E$ denotes completed homology $\pi_*L_{K(n)}(E \sma
E)$. Which $n$ is intended varies, but is generally clear from context.

$\G^u$ is the universal deformation formal group of $\Gamma$, defined over
$E_0$. In particular, there is a canonical isomorphism $\G^u \otimes_{E_0} k
\cong \Gamma$. We write $\H$ for the height $n-1$ formal group $\G^u \otimes
k((u_{n-1}))$. $\G$ will generally be used for other formal groups. If $f:R_1
\to R_2$ is a ring map and $\G$ is a formal group over $R_1$, then $f^*\G$ is
its base change over $R_2$ -- the upper star because this is the pullback of $\G
\to \Spec R_1$ along $\Spec R_2 \to \Spec R_1$.

We will write $\Def_\Gamma$ for the deformations functor of $\Gamma$, and
$\Def_\H^\aug$ for the functor of deformations of $\H$ augmented with
$\CLaur$-algebra structure; both are defined more explicitly below.

In section 6, $\psi^p$ and $\theta$ are certain operations on $\theta$-algebras.

\section{Background on deformation theory}\label{background}

In this introductory section, we review some facts from deformation theory, in
particular the Lubin-Tate theorem on deformations of formal groups.

\subsection{Witt vectors and Cohen rings}\label{witt vectors and cohen rings}

We will begin by presenting a deformation-theoretic point of view on the Witt
vectors of fields.

\begin{definition}\label{CLN}
  Let $\CLN$ be the category of complete noetherian local rings and continuous
  maps, where a ring $R$ with maximal ideal $\mathfrak{m}$ is equipped with its
  $\mathfrak{m}$-adic topology. More generally, if $A$ is a ring equipped with a
  chosen maximal ideal $I$, let $\CLN_A$ be the category whose objects are
  complete noetherian local $A$-algebras $R$ such that the structure map $A \to
  R$ sends $I$ into the maximal ideal of $R$, and whose morphisms are continuous
  $A$-algebra maps.
\end{definition}

We will not define the Witt vectors, but urge the reader to consult
\cite{Hesselholt}, \cite{Rabinoff}, or \cite{Serre} for more information. Let's
recall that if $R$ is a ring, its Witt vectors $W(R)$ are another ring equipped
with:
\begin{itemize}
\item a \textbf{Frobenius} map $F:W(R) \to W(R)$, which is a ring homomorphism;
\item a \textbf{Verschiebung} map $V:W(R) \to W(R)$, which is additive and
  satisfies $V(xF(y)) = V(x)y$ and $FV = p$;
\item a \textbf{multiplicative lift}, which is a group homomorphism $R^\times
  \to W(R)^\times$ written $x \mapsto [x]$.
\end{itemize}

We are mainly concerned with the Witt vectors of perfect fields, which are
described by the following statement.

\begin{proposition}[cf.\ \cite{Rabinoff}] If $k$ is a perfect field of
  characteristic $p$, then $W(k)$ is a $p$-torsion-free complete local ring with
  maximal ideal $VW(k) = pW(k)$ and residue field $k$. Every element of $W(k)$
  has a unique expansion of the form $\sum_{n=0}^\infty p^n[b_n]$, with $b_n \in
  R$.
\end{proposition}

The Witt vectors of a perfect characteristic $p$ field enjoy a universal
property, which we can describe as follows.

\begin{thm}\label{Wittuniv}
  Let $k$ be a perfect field of characteristic $p$ and let $R$ be a ring which
  is complete with respect to an ideal $I$ that contains $p$. Then for any map
  $i:k \to R/I$, there exists a unique continuous map completing any diagram of
  the form
  \[ \xymatrix{ Wk \ar@{-->}[r] \ar[d] & R \ar[d] \\ k \ar[r]_i & R/I. } \]
\end{thm}

\begin{proof}
  This proof originally goes back to Cartier, and one should consult \cite[\S
  II.4-6]{Serre}. Recall that elements of $Wk$ for $k$ perfect can be uniquely
  written in the form $\sum [a_n]p^n$, where $[a_n]$ is the Teichm\"uller lift
  of $a_n \in k$. The idea is that there is a unique multiplicative lift
  $\tau:k^\times \to R^\times$. One is then forced to send $\sum [a_n]p^n$ to
  $\sum \tau(a_n) p^n$, which converges by completeness of $R$.

  Regard $k$ as a subring of $R/I$. For each $a \in k$, define
  \[ U_n(a) = \{x^{p^n}:x \in R,\,x\equiv a^{p^{-n}}\pmod{I}\}. \] Here
  $a^{p^{-n}}$ is the unique $p^n$th root of $a$ in $k$. We have $U_{n+1}(a)
  \subseteq U_n(a)$. Moreover, if $x^{p^n}$ and $y^{p^n}$ are elements in
  $U_n(a)$, then $x \equiv y$ mod $I$, and thus $x^{p^n} \equiv y^{p^n}$ mod
  $I^{n+1}$ using the binomial theorem and the fact that $p \in I$. By
  completeness of $R$, there is a unique element in $\bigcap_{n \ge 0} U_n(a)$.
  Call this $\tau(a)$.

  One now observes that $\tau(a^p) = \tau(a)^p$, and that $\tau$ is the unique
  section $k^\times \to R^\times$ with this property. Indeed, if $\tau'$ also
  has this property, then
  \[ \tau'(a) = \tau'(a^{p^{-n}}) \in U_n(a)\text{ for all }n, \] so $\tau'(a) =
  \tau(a)$. Thus, there is \textit{at most one} multiplicative section. But
  $\tau$ is also multiplicative, because $U_n(a)\cdot U_n(b) \subseteq U_n(ab)$.
\end{proof}

\begin{remark}
  Taking $I$ to be the maximal ideal $\mathfrak{m}$ of a complete local
  noetherian ring $R$, one can see this as a trivial case of deformation theory.
  Indeed, define the functor of \textbf{deformations of nothing}, $\Def:\CLN \to
  \Sets$, by
  \[ \Def(R) = \{i:k \to R/\mathfrak{m}\}. \] Then we have proved that $\Def(R)
  \cong \Hom_{cts}(Wk, R)$.

  However, the above theorem proves slightly more. A map $k \to R/I$ gives
  deformations of nothing over numerous complete local rings, namely the
  completions of the localizations of $R$ at maximal ideals containing $I$. The
  theorem implies that these deformations of nothing assemble over $\Spf R$ to
  give a unique deformation of nothing over everything, which the reader must
  admit is really something.
\end{remark}

If $k$ is a non-perfect characteristic $p$ field, its ring of Witt vectors is
worse-behaved. There is still a surjection $\mathbf{w}_0: W(k) \to k$ with
kernel $VW(k)$, but this ideal need not be principal, and $V(1) \ne p$. In
addition, the universal property of \Cref{Wittuniv} is not satisfied. Instead,
one can construct rings having a weak version of this universal property.
(Topologists wishing to know more should also consult the last section of
\cite{AMS}.)

\begin{definition}
  A \textbf{Cohen ring} for a characteristic $p$ field $k$ is a complete
  discrete valuation ring with residue field $k$ and uniformizer $p$.
\end{definition}

\begin{example}
  The Witt vectors of a perfect field $k$ are a Cohen ring for $k$. For an
  imperfect field, we have $px = V(1)x = V(F(x))$, so the set of multiples of
  $p$ is in general a proper subset of the maximal ideal $VWk$. Thus, the Witt
  vectors are not a Cohen ring in this case.
\end{example}

\begin{example}
  If $k$ is perfect, the completed Laurent series ring $\CLaur =
  Wk((x))^\wedge_p$ is a Cohen ring for the Laurent series field $k((x))$.
\end{example}

\begin{thm}[{\cite[Tag 0323]{Stacks}}]\label{existence of cohen rings} Every
  characteristic $p$ field $k$ has a Cohen ring. If $C$ is a Cohen ring for the
  characteristic $p$ field $k$, then for every $n$, $\Z/p^n\Z \to C/p^nC$ is
  formally smooth.
\end{thm}

\begin{corollary}\label{existence of maps from cohen rings}
  If $\Coh$ is a Cohen ring for $k$, $R$ is a ring which is complete with
  respect to an ideal $I$ containing $p$, and $i:k \to R/I$ is an inclusion,
  then there exists a map completing the diagram
  \[ \xymatrix{ \Coh \ar@{-->}[r] \ar[d] & R \ar[d] \\ k \ar[r]_i & R/I. } \]
\end{corollary}

\begin{proof}
  Starting with the map $\Coh/p = k \to R/I$, we use formal smoothness to
  inductively construct lifts
  \[ \xymatrix{ \Z/p^n\Z \ar[d] \ar[r] & R/I^n \ar[d] \\ \Coh/p^n\Coh
      \ar@{-->}[ur] \ar[r] & R/I^{n-1} } \] along the square-zero extensions
  $R/I^n \to R/I^{n-1}$. By completeness, these assemble to a map $\Coh \to R$
  with the desired property.
\end{proof}

\begin{corollary}
  Any two Cohen rings for $k$ are (non-uniquely) isomorphic.
\end{corollary}

\begin{proof}
  Let $\Coh_1$ and $\Coh_2$ be Cohen rings for $k$. By the previous corollary,
  there exists a map $f:\Coh_1 \to \Coh_2$ reducing to the identity on $k$. If
  $f(x) = 0$, then $x$ must be divisible by $p$, and writing $x = px_1$ and
  proceeding inductively, we see that $x$ is divisible by all powers of $p$, so
  is zero by completeness. Thus, $f$ is injective. If $y \in \Coh_2$, then there
  is an $x_0 \in \Coh_2$ such that $f(x_0) - y$ is divisible by $p$, so is equal
  to some $py_1$. Proceeding inductively again and using completeness, we see
  that $f$ is surjective.
\end{proof}

\subsection{Formal groups and their deformations}

\begin{definition}
  By a \textbf{formal group} over a ring $R$, we mean a connected
  one-dimensional smooth commutative formal group scheme over $R$, which admits
  a global coordinate. If $R$ has a topology (for example, if $R$ is a complete
  local ring), then the formal scheme structure on any formal group is taken to
  be compatible with the topology. Thus, the underlying formal scheme of any
  formal group over $R$ is non-canonically isomorphic to $\Spf R[[x]]$.

  Given a ring homomorphism $R \to R'$ (continuous, if the rings involved have
  topologies), and a formal group $\Gamma$ over $R$, we write $\mathbb{G}
  \otimes_R R'$ for the base change of $\Gamma$ to $R'$.
\end{definition}

Formal groups over a field of characteristic $p$ are stratified by an invariant
known as height, which can be briefly defined as follows. Choosing a coordinate
$x$ for the formal group $\Gamma$ and writing $[p]_\Gamma(x)$ for the power
series that expresses formal multiplication by $p$, one finds that
\[
  [p]_\Gamma(x) = ux^{p^n} + \dotsb
\]
for some $0 < n \le \infty$ independent of the choice of $x$. This $n$ is the
\textbf{height} of $\Gamma$. We will only be concerned with formal groups of
finite height (i.~e.~where $[p]_\Gamma(x)$ is not identically 0) in this paper.

We write $\Mfg$ for the moduli of formal groups, which is an algebraic stack.
For further background on formal groups, see \cite[Appendix 2]{Rav} and
\cite{Goerss}.

\begin{definition}
  Let $\Gamma$ be a formal group over $k$ and $R$ a complete local noetherian
  ring with maximal ideal $\mathfrak{m}$. A \textbf{deformation} of $\Gamma$
  over $R$ is a triple
  \[ (\mathbb{G}, \overline{i}, \alpha), \] where
  \begin{itemize}
  \item $\mathbb{G}$ is a formal group over $R$,
  \item $\overline{i}$ is an inclusion $k \to R/\mathfrak{m}$,
  \item and $\alpha$ is an isomorphism $\Gamma \otimes_k^{\overline{i}}
    R/\mathfrak{m} \stackrel{\sim}{\to} \mathbb{G} \otimes_R R/\mathfrak{m}$.
  \end{itemize}

  An \textbf{isomorphism} of deformations of $\Gamma$ over $R$,
  $\phi:(\mathbb{G}_1, \overline{i_1}, \alpha_1) \to (\mathbb{G}_2,
  \overline{i_2}, \alpha_2)$ is
  \begin{itemize}
  \item the condition that $\overline{i_1} = \overline{i_2}$,
  \item and an isomorphism $\phi:\mathbb{G}_1 \to \mathbb{G}_2$ of formal groups
    over $R$,
  \item such that the square
    \[ \xymatrix{ \Gamma \otimes_k R/\mathfrak{m} \ar[d]^1 \ar[r]^{\alpha_1} &
        \G_1 \otimes_R R/\mathfrak{m} \ar[d]_{\phi} \\ \Gamma \otimes_k
        R/\mathfrak{m} \ar[r]_{\alpha_2} & \G_2 \otimes_R R/\mathfrak{m} } \]
    commutes.
  \end{itemize}

  If $\Gamma$ is a formal group over $k$, then let
  \[ \Def_\Gamma:\CLN \to \Gpd \] be the functor that sends $R$ to the groupoid
  of deformations of $\Gamma$ over $R$ and their isomorphisms. More generally,
  if $A$ is a ring with a maximal ideal $I$ such that $A/I = k$, and $\Gamma$ is
  a formal group over $k$, then let
  \[ \Def_\Gamma^A:\CLN_A \to \Gpd \] be the functor that sends an $A$-algebra
  $R$ to the groupoid of deformations $(\G, \overline{i}, \alpha)$ of $\Gamma$
  over $R$ such that $\overline{i}:A/I \to R/\mathfrak{m}_R$ is the reduction of
  the $A$-algebra structure map, and their isomorphisms.
\end{definition}

The Lubin-Tate theorem is originally stated in the following form \cite{LT}:

\begin{thm}[Lubin-Tate]\label{Lubin-Tate} Let $\Gamma$ be a formal group of
  finite height $n$ over a field $k$ of characteristic $p > 0$, let $A$ be a
  ring with an ideal $I$ such that $A/I \cong k$, and suppose given a diagram
  \[
    \xymatrix{ A \ar[r]^i \ar[d] & R \ar[d] \\ k \ar[r]_{\overline{i}} &
      R/\mathfrak{m}}
  \]
  where $R$ is a complete local noetherian ring with maximal ideal
  $\mathfrak{m}$. Then there is a formal group $\mathbb{G}^u$ over $A[[u_1,
  \dotsc, u_{n-1}]]$ such that, for any formal group $\mathbb{G}$ over $R$ with
  an isomorphism
  \[ \alpha:\Gamma \otimes_k R/\mathfrak{m} \stackrel{\sim}{\to} \mathbb{G}
    \otimes_R R/\mathfrak{m}, \] there is a unique continuous $A$-algebra map
  $f:A[[u_1,\dotsc, u_{n-1}]] \to R$ and a unique isomorphism
  \[ \mathbb{G}^u \otimes_{A[[u_1,\dotsc,u_{n-1}]]} R \stackrel{\sim}{\to}
    \mathbb{G} \] that reduces to $\alpha$ over $R/\mathfrak{m}$.
\end{thm}

In the language defined above, we can restate this as follows.

\begin{corollary}[Lubin-Tate theorem, equivalent form] Let $\Gamma$ be a formal
  group of finite height $n$ over a field $k$ of characteristic $p > 0$. Then
  the functor $\Def_\Gamma^A$ on $\CLN_A$ is represented by $A[[u_1, \dotsc,
  u_{n-1}]]$. In other words, for $R \in \CLN_A$, the groupoid
  $\Def_\Gamma^A(R)$ is naturally equivalent to the discrete groupoid
  $\Maps_{A,cts}(A[[u_1,\dotsc,u_{n-1}]], R)$.
\end{corollary}

Typically, one is interested in representing the functor $\Def_\Gamma$ on the
category $\CLN$, not on the categories $\CLN_A$. Using \Cref{Wittuniv}, this
version of the theorem can be proved if the field $k$ is perfect (and most
often, this refinement is the one used by homotopy theorists).

\begin{corollary}[Lubin-Tate theorem for perfect fields]\label{Lubin-Tate
    perfect} Let $\Gamma$ be a formal group of finite height $n$ over a
  {\upshape perfect} field $k$ of characteristic $p > 0$. Then the functor
  $\Def_\Gamma$ on $\CLN$ is represented by $Wk[[u_1, \dotsc, u_{n-1}]]$.
\end{corollary}

\begin{proof}
  Given $(\G, \overline{i}, \alpha) \in \Def_\Gamma(R)$, let $i:Wk \to R$ be the
  unique continuous lift over $\overline{i}$ guaranteed by \Cref{Wittuniv}. The
  map $i$ makes $R$ an object of $\CLN_{Wk}$, and $(\G, \overline{i}, \alpha)$
  an object of $\Def_\Gamma^{Wk}(R)$. This proves that the inclusion map
  $\Def_\Gamma^{Wk}(R) \to \Def_\Gamma(R)$ is essentially surjective, but
  $\Def_\Gamma^{Wk}(R)$ is a full subgroupoid of $\Def_\Gamma(R)$ by definition.
  So the inclusion is an equivalence, and on the other hand we have
  \[ \Maps_{Wk,cts}(Wk[[u_1,\dotsc,u_{n-1}]],R) =
    \Maps_{cts}(Wk[[u_1,\dotsc,u_{n-1}]],R), \] again by \Cref{Wittuniv}.
\end{proof}

\begin{remark}\label{parameters}
  The ring $Wk[[u_1, \dotsc, u_{n-1}]]$ is called the \textbf{Lubin-Tate ring}
  for $(k, \Gamma)$. As a result of \Cref{Lubin-Tate perfect}, it carries a
  deformation $(\G^u, 1, \alpha^u)$ which is a universal deformation of
  $\Gamma$, in the sense that any other deformation $(\G, i, \alpha)$ over $R
  \in \CLN$ is uniquely isomorphic to the base change of $\G^u$ along a unique
  map $Wk[[u_1, \dotsc, u_{n-1}]] \to R$.

  In fact, this can be made fairly explicit. If $\Gamma$ has the Honda formal
  group law over $k$, with $p$-series
  \[ [p]_\Gamma(x) = x^{p^n} \] (for some chosen coordinate $x$), then we can
  choose a coordinate on $\G^u$ such that
  \[ [p]_{\G^u}(x) = px +_{\G^u} u_1x^p +_{\G^u} \dotsb +_{\G^u}
    u_{n-1}x^{p^{n-1}} +_{\G^u} x^{p^n}. \] The isomorphism $\alpha^u:\Gamma
  \stackrel{\sim}{\to} \G^u \otimes k$ matches these two coordinates.

  In other words, a deformation of a height $n$ formal group is specified, up to
  isomorphism, by deformations of the coefficients of $x^p, \dotsc, x^{p^{n-1}}$
  in its $p$-series, and the Lubin-Tate parameters keep track of these
  deformations. This should indicate that the parameters $u_1, \dotsc, u_{n-1}$
  themselves are very non-canonical, and are related to the choice of coordinate
  on the formal group $\Gamma$.
\end{remark}

\subsection{Automorphisms of formal groups}

\begin{definition}\label{morava stabilizer}
  Let $\Gamma$ be a height $n$ formal group over a field $k$. The \textbf{Morava
    stabilizer group} $\Aut(k,\Gamma)$ is the group of pairs $(\tau, g)$, where
  $\tau:k \stackrel{\sim}{\to} k$ is an automorphism, and $g$ is an isomorphism
  of formal groups over $k$, $g:\Gamma \stackrel{\sim}{\to} \tau^*\Gamma$.
\end{definition}

Equivalently, one can write $(\tau, g) \in \Aut(k,\Gamma)$ as a commutative
square
\[ \xymatrix{ \Gamma \ar[r]^{\sim} \ar[d] & \Gamma \ar[d] \\ \Spec k
    \ar[r]_{\tau^*} & \Spec k } \] where the horizontal maps are isomorphisms,
the top one a (non-$k$-linear) isomorphism of formal groups. Then $g:\Gamma
\stackrel{\sim}{\to} \tau^*\Gamma$ is the map induced by the universal property
of the pullback.

The composition law can be written
\[ (\tau_2, g_2)(\tau_1, g_1) = (\tau_2\tau_1, \tau_2^*(g_1)g_2). \] The odd
variance of this formula is a result of our choice to write $\tau$ as a map of
fields (that is, rings) and $g$ as a map of formal groups (that is, formal
schemes).

\begin{remark}\label{morava stabilizer action on Def}
  There is a left action of $\Aut(k,\Gamma)$ on $\Def_\Gamma$, defined as
  follows. Given $(\G, i, \alpha) \in \Def_\Gamma(R)$ and $(\tau:k \to k,
  g:\Gamma \to \tau^*\Gamma) \in \Aut(k,\Gamma)$ (see \Cref{morava stabilizer}),
  define
  \[ (\tau,g)(\G, i, \alpha) = (\G, i \circ \tau, \alpha g^{-1}:\Gamma
    \otimes_k^{i\tau} R/\mathfrak{m} \stackrel{g^{-1}}{\to} \Gamma \otimes_k^i
    R/\mathfrak{m} \stackrel{\alpha}{\to} \G \otimes_R R/\mathfrak{m}). \]
\end{remark}

\begin{remark}
  Typically, in homotopy theory, $k$ is taken to be a finite field containing
  $\F_{p^n}$. Then all automorphisms of $\Gamma$ are defined over $k$, i.~e.,
  $\Aut_k(\Gamma) = \Aut_{\overline{\F_p}}(\Gamma)$. It follows that
  \[ \Aut(k,\Gamma) = \Gal(k/\F_p) \ltimes \Aut_k(\Gamma). \] Moreover,
  $\Aut_k(\Gamma)$ is a profinite group, isomorphic to the group of units in a
  maximal order of the division algebra of invariant $1/n$ over $Wk[1/p]$. One
  commonly writes $\mathbb{S}_n$ for $\Aut_k(\Gamma)$ and $\mathbb{G}_n$ for
  $\Aut(k, \Gamma)$ in this case.
\end{remark}

Now that we have described the automorphisms of a formal group over a field, we
describe how they deform. While deformations of formal groups are controlled by
Lubin-Tate parameters, isomorphisms of formal groups deform uniquely.

\begin{definition}
  Let $\G_1$ and $\G_2$ be two formal groups over rings a ring $A$. The
  \textbf{moduli of isomorphisms} from $\G_1$ to $\G_2$ is the functor
  $\Iso(\G_1, \G_2):\Alg_A \to \Sets$ given by
  \[ \Iso(\G_1, \G_2)(B) = \{\phi:\G_1 \otimes_{A} B \stackrel{\sim}{\to} \G_2
    \otimes_{A} B\text{ an isomorphism of formal groups }\}. \]
\end{definition}

\begin{thm}[{\cite[5.25]{Goerss}}]\label{etaleness of isomorphisms} If
  $\Gamma_1$ and $\Gamma_2$ are formal groups over a ring $A$, then $\Iso(\G_1,
  \G_2)$ is an affine $A$-scheme. Moreover, if $\G_1$ and $\G_2$ are height $n$
  formal groups over a field $k$ of characteristic $p$, then $\Iso(\G_1, \G_2)$
  is pro-\'etale over $k$.
\end{thm}

We will use the following easy corollary repeatedly.

\begin{corollary}
  Let $\G_1$ and $\G_2$ be height $n$ formal groups over a field $k$ of
  characteristic $p$, and let $R$ be a a $k$-algebra which is complete with
  respect to an ideal $I$. Then $\Iso_R(\G_1, \G_2) \to \Iso_{R/I}(\G_1, \G_2)$
  is an isomorphism.
\end{corollary}

\begin{proof}
  If $I$ is nilpotent, this is just the infinitesimal criterion of \'etaleness.
  In general, we write $R = \lim R/I^n$ and note that the moduli of isomorphisms
  commutes with limits.
\end{proof}

\section{Background on homotopy theory}

\subsection{Morava \texorpdfstring{$E$}{E}-theory}

The Lubin-Tate rings can be realized in homotopy theory in an interesting way,
as the coefficient rings of cohomology theories known as \textbf{Morava
  $E$-theory}.

\begin{thm}
  Let $\Gamma$ be a height $n$ formal group over a perfect field $k$ of
  characteristic $p$. There is a complex orientable, even periodic ring spectrum
  $E = E(k,\Gamma)$ such that $\pi_0E = Wk[[u_1, \dotsc, u_{n-1}]]$ and the
  formal group of $E$ (rescaled to degree zero) is the universal deformation of
  $\Gamma$.
\end{thm}

One way to prove this is using the Landweber exact functor theorem, which
applies to the Lubin-Tate rings because they are flat over the moduli of formal
groups \cite{Goerss}.

There are also more homotopical constructions available. For instance, the
Brown-Peterson spectrum $BP$ is a summand of $p$-local complex cobordism, with
coefficient ring
\[ BP_* = \Z_{(p)}[v_1, v_2, \dotsc], \quad |v_i| = 2(p^i-1). \] One can define
a ring spectrum $BP\langle n\rangle$ as the quotient of $BP$ by the ideal
$(v_{n+1}, v_{n+2}, \dotsc)$ \cite{EKMM}. The Johnson-Wilson $E$-theory $E(n)$
is then $BP\langle n\rangle[v_n^{-1}]$. This admits an \'etale extension in ring
spectra, $\widetilde{E(n)}$, with
\[ \pi_*\widetilde{E(n)} = \pi_*E(n)[u]/(u^{(p^n-1)} - v_n) =
  \Z_{(p)}[u_1,\dotsc, u_{n-1}][u^{\pm 1}] \] where $|u| = 2$ and $|u_i| = 0$,
and $u_i = v_iu^{1-p^i}$. Finally, the $K(n)$-localization of $\widetilde{E(n)}$
is a Morava $E$-theory. Specifically, it is the Morava $E$-theory for the Honda
formal group of height $n$ over $\F_p$, with
\[ [p]_\Gamma(x) = x^{p^n}. \]

Both constructions lack a certain something. Both the Landweber exact functor
theorem, and the process of quotienting ring spectra by ideals, only produce
Morava $E$-theory as a homotopy commutative ring spectrum. More advanced
obstruction theory techniques imply the following, for which one should see
\cite{GH1} and \cite{RezkHM}.

\begin{thm}[Goerss-Hopkins-Miller] There is a unique $\Einf$ ring spectrum, up
  to $\Einf$ homotopy equivalence, whose underlying ring spectrum is
  $E(k,\Gamma)$. Moreover, the space of $\Einf$ endomorphisms of this spectrum
  is homotopy equivalent to the discrete group $\Aut(k,\Gamma)$.
\end{thm}

In other words, by insisting on more structure -- the structure of an $\Einf$
ring spectrum -- we are able to realize the Lubin-Tate formal groups in stable
homotopy theory in an essentially unique way. Thus, there is an unexpectedly
tight correspondence between chromatic stable homotopy theory and the algebraic
geometry of formal groups, so long as one works locally on both sides: Morava
$E$-theories are $K(n)$-local, while the Lubin-Tate formal groups are universal
deformations of formal groups over fields. In turn, this suggests that, as a
general principle, \textit{all} phenomena having to do with Morava $E$-theory
can be described in terms of formal groups.

\subsection{L-completeness and completed \texorpdfstring{$E$}{E}-homology}

The $E$-theory of a point -- the Lubin-Tate ring -- is a complete local ring.
The $E$-theory spectrum itself also satisfies a completeness property, in that
it is $K(n)$-local: this implies that, letting $I = (i_0, \dotsc, i_{n-1})$
range over a cofinal sequence of $n$-tuples of integers such that the
generalized Moore spectrum $S/(p^{i_0}, v_1^{i_1}, \dotsc, v_{n-1}^{i_{n-1}})$
is defined, the natural map
\[ E \to \holim_I E \sma S/(p^{i_0}, \dotsc, v_{n-1}^{i_{n-1}}) \] is an
equivalence. One expects to be able to deduce completeness properties for the
$E$-homology and $E$-cohomology functors as well. For example, for any spectrum
$X$,
\begin{multline*}
  E^*X = \pi_*F(X,E) = \pi_*\holim_I F(X, E \sma S/(p^{i_0}, \dotsc,
  v_{n-1}^{i_{n-1}})) \\
  = \pi_*\holim_I F(X, E) \sma S/(p^{i_0}, \dotsc, v_{n-1}^{i_{n-1}}).
\end{multline*}
This is an $E_*$-module with a derived completeness property to be explained
momentarily.

The $E$-homology groups $E_*X = \pi_*(E \sma X)$ do not generally have this
completeness property, because the smash product does not distribute over the
homotopy limit. Instead, we work with \textbf{completed $E$-homology},
\[
  E_*^\wedge X = \pi_*L_{K(n)}(E \sma X).
\]
This is not a homology theory in the Eilenberg-Steenrod sense. Nevertheless, it
is often better behaved than uncompleted $E$-homology: for example, while $E_*E$
is quite complicated, $E_*^\wedge E$ is pro-free over $E$ and has the simple
Hopf algebroid description $\Maps_{cts}(\G_n, E_*)$.

We now describe the derived completeness property of $E$-cohomology and
completed $E$-homology, known as L-completeness. More information, and proofs of
the results mentioned here, can be found in \cite{GM}, \cite[Appendix A]{HS},
and \cite{BF}. In practice, the distinction between L-completeness and classical
completeness will not matter much here, as any L-complete module we will ever
consider will be classically complete.

\begin{definition}
  Let $R$ be a Noetherian ring and let $I$ be an ideal in $R$, generated by a
  regular sequence of length $n$. We write $M^\wedge$ for the $I$-adic
  completion of an $R$-module $M$, and $L_sM$ for the $s$th derived functor of
  $I$-adic completion applied to $M$. There are natural maps
  \[ M \to L_0 M \to M^\wedge. \] A module $M$ is \textbf{L-complete} if the
  natural map $M \to L_0M$ is an isomorphism, and \textbf{(classically)
    complete} if the map $M \to M^\wedge$ is an isomorphism.

  Write $\CplMod_R$ for the category of L-complete $R$-modules.
\end{definition}

In practice, we will usually take $R = E_*$, $I = I_n = (p, u_1, \dotsc,
u_{n-1})$, and work with the degreewise versions of these properties for graded
rings. The hypotheses on $R$ are less general than those in \cite{GM} (which
weakens the Noetherian condition) and more general than those in \cite{HS} and
\cite{BF} (which require $R$ to be a local ring and $I$ its maximal ideal).
However, the proofs in \cite{HS} and \cite{BF} do not rely on $I$ being maximal.

The following proposition summarizes results from \cite{GM} and \cite{HS}.

\begin{proposition}
  \begin{enumerate}
  \item $\CplMod_R$ is an exact subcategory of $\Mod_R$, which is closed under
    extensions, limits, and $\lim^1$ of towers.
  \item There is an adjunction $L_0:\Mod_R \leftrightarrows \CplMod_R:i$, where
    $i$ is the inclusion.
  \item If $M$ is free, or finitely generated, then the natural map $M^\wedge
    \to L_0M$ is an isomorphism.
  \item If $M$ is classically complete, it is L-complete.
  \item In general, $L_sM = 0$ for $s \ge n+1$, and there are natural exact
    sequences
    \[ 0 \to \lim_k\nolimits^1 \Tor^R_{s+1}(R/I^k, M) \to L_s M \to \lim_k
      \Tor^R_s(R/I^k, M) \to 0. \]
  \end{enumerate}
\end{proposition}

\begin{proposition}[{\cite[Corollary 2.3, Proposition 8.4]{HS}}] The functors
  $E^*$ and $E^\wedge_*$ are naturally valued in the category $\CplMod_{E_*}$ of
  graded $E_*$-modules which are L-complete with respect to the ideal $I_n = (p,
  u_1, \dotsc, u_{n-1})$.
\end{proposition}

\begin{proposition}[{\cite[Proposition 8.4]{HS}}] If $X$ is finite, then
  $E_*^\wedge X = E_*X$. If $E_*X$ is $E_*$-free, then $E_*^\wedge X$ is its
  $I_n$-adic completion.
\end{proposition}

\begin{proof}
  We repeat the proof from \cite{HS} because it will be relevant later on. If
  $X$ is finite, then $E \sma X$ is an iterated cofiber of suspensions of $E$,
  so is already $K(n)$-local. For general $X$, $v_n$ is invertible on $E \sma
  X$, so we have
  \[ L_{K(n)}(E \sma X) = \holim_I E \sma X \sma S/(p^{i_0}, \dotsc,
    v_{n-1}^{i_{n-1}}). \] If $E_*X$ is $E_*$-flat, then
  \[ \pi_*(E \sma X \sma S/(p^{i_0}, \dotsc, v_{n-1}^{i_{n-1}})) =
    E_*X/(p^{i_0}, \dotsc, v_{n-1}^{i_{n-1}}), \] and the transition maps in the
  inverse system of homotopy groups are surjective. Therefore,
  \[ \pi_*(\holim_I E \sma X \sma S/(p^{i_0}, \dotsc, v_{n-1}^{i_{n-1}})) =
    \lim_I E_*X/(p^{i_0}, \dotsc, v_{n-1}^{i_{n-1}}) = (E_X)^\wedge_{I_n}. \]
\end{proof}

\subsection{Cooperations for \texorpdfstring{$E$}{E}-theory}

The completed $E$-homology $E^\wedge_*X = \pi_*L_{K(n)}(E \sma X)$ of a space or
spectrum $X$ is naturally a complete comodule for a coalgebra of cooperations
$E^\wedge_*X$. This has a surprisingly simple form. In this section, we write
$\G_n = \Aut(k,\Gamma)$.

\begin{thm}[{\cite{DH}}]\label{E cooperations} If $E = E(k,\Gamma)$ where
  $\Gamma$ is the height $n$ Honda formal group over a perfect field $k$ of
  characteristic $p$, then there is a Hopf algebroid isomorphism
  \[ E^\wedge_* E \cong \Cts(\G_n, E_*), \] where $\Cts(\G_n, E_*)$ is the
  $E_*$-algebra of continuous set maps $\G_n \to E_*$.
\end{thm}

As an immediate application, if $X$ is a spectrum such that $E_*^\wedge X$ is
honestly $I_n$-complete, then the $E_*^\wedge E$-comodule structure on
$E_*^\wedge X$ is equivalently given by a continuous $\G_n$ action. (If
$E_*^\wedge X$ is merely L-complete, then a comodule structure should instead be
taken as the \textit{definition} of ``continuous action''.) Thus in the
$K(n)$-local $E$-based Adams spectral sequence
\[ E_2 = \Ext_{E^\wedge_* E}^*(E_*, E^\wedge_*X) \Rightarrow \pi_*L_{K(n)}X, \]
we can rewrite the $E_2$ page as group cohomology $H^*_{cts}(\G_n,
E^\wedge_*X)$.

For the sake of inspiring the arguments below, we give a proof of this theorem.
Let $\Alg_{E_*}^\wedge$ be the category of even periodic, $I_n$-adically
complete $E_*$-algebras, and let $\Alg_{E_0}^\wedge$ be the category of
$I_n$-adically complete $E_0$-algebras. The proof will show that $E_*^\wedge E$
and $\Cts(\G_n, E_*)$ represent the same functor on $\Alg_{E_*}^\wedge$.

\begin{lemma}\label{E cooperations periodic}
  The ring $E_*^\wedge E$ is an object in $\Alg_{E_*}^\wedge$.
\end{lemma}

\begin{proof}
  Since $E$ is Landweber exact,
  \[ E_*E = \pi_*(E \sma E) = E_* \otimes_{BP_*} BP_*BP \otimes_{BP_*} E_* \] is
  even periodic and flat as a left $E_*$-module. The theory $E \sma E$ is
  $L_n$-local, which means that
  \[ L_{K(n)}(E \sma E) = \holim_I E \sma E \sma S/I = \holim_I E \sma E/I \]
  where $I$ ranges over ideals $(p^{i_0}, v_1^{i_1}, \dotsc, v_{n-1}^{i_{n-1}})$
  such that the associated generalized Moore spectrum exists. Again, the
  homotopy groups of the objects in the limit diagram are
  \[ E_*(E/I) = E_* \otimes_{BP_*} BP_*BP \otimes_{BP_*} E_*/I. \] Thus, the
  maps in the diagram are surjective on homotopy groups, so that
  \[ \pi_*L_{K(n)}(E \sma E) = \lim_I E_*(E/I) = (E_* \otimes_{BP_*} BP_*BP
    \otimes_{BP_*} E_*)^\wedge_{I_n}. \] As $I_n$ and its powers are images of
  invariant ideals in $BP_*BP$, it doesn't matter whether we complete with
  respect to the $I_n$ coming from the left or right $E_*$-module structure.
  Clearly, this is an object of $\Alg_{E_*}^\wedge$.
\end{proof}

\begin{lemma}\label{E cooperations continuity}
  For $R_* \in \Alg_{E_*}^\wedge$, pre-composition with the completion map $E_*E
  \to E_*^\wedge E$ induces an isomorphism
  \[ \Hom_{\Alg_{E_*}^\wedge}(E_*^\wedge E, R_*) \cong \Hom_{E_*}(E_*E, R_*). \]
  Moreover, for any map $E_*^\wedge E \to R_*$ in $\Alg_{E_*}^\wedge$, the
  composition
  \[ E_* \stackrel{\eta_R}{\to} E_*E \to R_* \] is also continuous.
\end{lemma}

\begin{proof}
  Let $R_* \in \Alg_{E_*}^\wedge$ and give $R_*$ the $I_n$-adic topology. As we
  saw in the proof of the previous lemma, $I_n$ is the image of an invariant
  ideal in $BP_*BP$. Thus, any map
  \[ f:E_* \otimes_{BP_*} BP_*BP \otimes_{BP_*} E_* \to R_* \] extending the
  given map $E_* \to R_*$ has $f(I_n\cdot E_*E) \subseteq I_n R_*$. In
  particular, $R_*$ is also complete with respect to $I_n\cdot E_*E$, so that
  $f$ factors uniquely through the completion $E^\wedge_* E$. Moreover,
  $f(\eta_R(I_n\cdot E_*))$ is also in $I_n R_*$, so that the map $E_* \to R_*$
  coming from the right unit is also continuous.
\end{proof}

Let $\Alg_{E_*,loc}^\wedge$ be the full subcategory of $R \in \Alg_{E_*}^\wedge$
such that $R_0$ is complete local. Then the map $E_0 \to R_0$ classifies an
object $(\G, i, \alpha) \in \Def_\Gamma(R_0)$.

\begin{proposition}\label{E cooperations representability}
  Let $R_* \in \Alg_{E_*,loc}^\wedge$, and let $(\G, i, \alpha)$ be the
  deformation of $\Gamma$ classified by $E_0 \to R_0$. Then the set of maps
  $\Hom_{\Alg_{E_*}^\wedge}(E^\wedge_*E, R_*)$ is naturally isomorphic to the
  set of pairs $(j,\gamma)$, where $j$ is a map $k \to R_0/\mathfrak{m}$, and
  $\gamma$ is an isomorphism of formal groups over $R/\mathfrak{m}$,
  $\gamma:\Gamma \otimes^i_k R/\mathfrak{m} \to \Gamma \otimes^j_k
  R/\mathfrak{m}$.
\end{proposition}

\begin{proof}
  By \Cref{E cooperations continuity},
  \[ \Hom_{\Alg_{E_*}^\wedge}(E_*^\wedge E, R_*) \cong \Hom_{E_*}(E_*E, R_*)
    \cong \Hom_{E_0}(E_0E, R_0). \]

  Write $BPP$ for 2-periodic $BP$. Then we also have
  \[ E_0E = E_0 \otimes_{BPP_0} BPP_0BPP \otimes_{BPP_0} E_0. \] The Hopf
  algebroid $(BPP_0, BPP_0BPP)$ presents the moduli of $p$-local formal groups,
  so there is a pullback square of stacks
  \[ \xymatrix{ \Spec E_0E \ar[r] \ar[d] \pullbackcorner & \Spec E_0 \ar[d] \\
      \Spec E_0 \ar[r] & \Mfg. } \] Again, \Cref{E cooperations continuity}
  implies that, for any $E_0$-map $E_0E \to R_0$, the map $E_0 \to R_0$ coming
  from the right unit is also continuous, and thus classifies a deformation. So
  we have a homotopy pullback of groupoids,
  \[ \xymatrix{ \Hom_{\Alg_{E_*}^\wedge}(E_*^\wedge E, R_*) \ar[r] \ar[d]
      \pullbackcorner & \Def_\Gamma(R_0) \ar[d] \\ \{*\} \ar[r]_{\G} &
      \Mfg(R_0). } \] An object in the pullback is given by another deformation
  $(\G', j, \beta) \in \Def_\Gamma(R_0)$ and an isomorphism $\phi:\G \to \G'$.
  An isomorphism in the pullback is an isomorphism $\psi:(\G', j, \beta) \to
  (\G'', j, \delta)$ such that the obvious triangle involving $\phi$ commutes.
  Now, there is an isomorphism of deformations $\phi^{-1}:(\G', j, \beta) \to
  (\G, j, \beta \circ \phi^{-1})$; this is the only isomorphism from $(\G', j,
  \beta)$ to a deformation whose underlying formal group is exactly $\G$. It
  follows that the connected components of the pullback groupoid are
  contractible, as expected, and correspond to pairs
  \[ (j, \beta:\Gamma \otimes^j R_0/\mathfrak{m} \stackrel{\sim}{\to} \G \otimes
    R_0/\mathfrak{m}). \] Equivalently, they correspond to pairs
  \[ (j:k \to R_0/\mathfrak{m}, \gamma = \beta^{-1}\alpha:\Gamma \otimes^i
    R_0/\mathfrak{m} \stackrel{\sim}{\to} \Gamma \otimes^j R_0/\mathfrak{m}). \]
\end{proof}

\begin{example}
  The group $\Aut(k,\Gamma)$ acts on this set by pre-composing with the map $j$
  and post-composing with the isomorphism $\gamma$. However, this action need
  not be transitive. For example, let $\Gamma$ be the height 1 formal group over
  the perfect field $K = \overline{\F_p}((u^{1/p^\infty}))$ with $p$-series
  \[ [p]_\Gamma(x) = ux^p. \] Let $L = K \otimes_{\overline{\F_p}} K =
  \overline{\F_p}((u^{1/p^\infty}, v^{1/p^\infty}))$, and let $R_0$ be the
  algebraic closure of $L$. There are two maps $j_1, j_2:K \to R_0$,
  respectively sending $u$ to $u$ and to $v$. The base changes of $\Gamma$ along
  these maps are isomorphic over $R_0$ via
  \[ \phi(x) = (u/v)^{1/(p-1)}x. \] This isomorphism is not induced by an
  element of $\Aut(k,\Gamma)$.
\end{example}

We now specialize to the case where $\Gamma$ is the height $n$ Honda formal
group over a finite field $k$ containing $\F_{p^n}$, with $[p]_\Gamma(x) =
x^{p^n}$. In this case, the formal group is algebraic enough to prevent the
above subtlety from occurring. (In fact, the following argument works in
slightly more generality: one can take $[p]_\Gamma(x) = ux^{p^n}$ for $u \in
k^\times$, which at least implies that $\Frob_\Gamma$ is central in
$\End_k(\Gamma)$.) The author thanks Paul Goerss for pointing out this subtlety
and the following method of addressing it.

\begin{proof}[Proof of \Cref{E cooperations}] First, we need to construct a
  continuous $E_*$-algebra map $f:E_*^\wedge E \to \Cts(\G_n, E_*)$. Such a map
  is adjoint to a continuous map
  \[ \G_n \to \Hom_{\Alg_{E_*}^\wedge}(E_*^\wedge E, E_*). \] Let $(\G^u, 1,
  \alpha^u)$ be the universal deformation over $E_*$. The $E_*$-algebra
  structure map $E_* \to E_*$ is just the identity map, which classifies this
  deformation. By \Cref{E cooperations representability},
  \[ \Hom_{\Alg_{E_*}^\wedge}(E_*^\wedge E, E_*) \cong \{(j:k \to k,
    \gamma:\Gamma \stackrel{\sim}{\to} \Gamma \otimes^j k)\}. \] Since $k$ is a
  finite field, any map $k \to k$ is an isomorphism. So the right-hand side is
  exactly $\G_n$, defining the desired map.

  For simplicity's sake, we now restrict everything to degree zero. To show that
  the map $f:E_0^\wedge E \to \Cts(\G_n, E_0)$ is an isomorphism, it suffices,
  since both sides are flat and complete $E_0$-algebras, that it induces an
  isomorphism mod $I_n$. Now, $I_n$ is an invariant ideal in $BP$, so
  \[ E_0^\wedge E/I_n = k \otimes_{BPP_0} BPP_0BPP \otimes_{BPP_0} k. \] A map
  from this into a ring $R$ is the same as a pair of maps $i,j:k \to R$ and an
  isomorphism $\gamma:\Gamma \otimes^i_k R \to \Gamma \otimes^j_k R$ of formal
  groups over $R$. Now, if $\Gamma$ is the Honda formal group over $k \supseteq
  \F_{p^n}$, we have a coordinate $x$ for $\Gamma$ with $[p]_\Gamma(x) =
  x^{p^n}$, and this must commute in the obvious way with any isomorphism
  $\gamma$. It follows that the coefficients of $\gamma$, viewed as a power
  series in $x$, are fixed by the $n$th power of Frobenius. Since $R$ is an
  $\F_{p^n}$-algebra via $i$, the subring of elements of $R$ fixed by the $n$th
  power of Frobenius is a product of copies of $\F_{p^n}$. By \Cref{etaleness of
    isomorphisms}, the isomorphism $\gamma$ is defined over $\F_{p^n}$.

  Thus, the data $(i, j, \gamma)$ is always base changed from data of the form
  \[ (1:k \to k, j:k \to k, \gamma: \Gamma \stackrel{\sim}{\to} j^*\Gamma). \]
  Since $k$ is finite, $j$ is an isomorphism. This is precisely an element of
  the Morava stabilizer group $\G_n = \Aut(k, \Gamma)$ (cf.\ \Cref{morava
    stabilizer}). Thus, the map $E_0^\wedge E/I_n \to \Cts(\G_n, k)$ is an
  isomorphism.
\end{proof}

\section{Localized \texorpdfstring{$E$}{E}-theory and augmented
  deformations}\label{LE}

Our primary concern in this paper is the spectrum $L_{K(n-1)}E_n$, where $n \ge
2$. We abbreviate this spectrum by $\LE$.

\subsection{Coefficients}
We begin by describing the coefficient ring of $\LE$.

\begin{proposition}\label{coefficients}
  The coefficient ring $\LE_*$ is even periodic, with
  \[ \LE_0 = Wk[[u_1, \dotsc, u_{n-1}]][u_{n-1}^{\pm 1}]^{\wedge}_{(p, u_1,
      \dotsc, u_{n-2})}. \]
\end{proposition}

\begin{proof}
  By \cite{Rav84}, $BP$ satisfies the telescope conjecture, in the sense that
  there is an equality of Bousfield classes $\langle BP \wedge T(n-1)\rangle =
  \langle BP \wedge K(n-1)\rangle$, where $T(n)$ is a $v_{n-1}$-telescope of a
  finite type $n-1$ spectrum. As $E_n$ is a $BP$-module, it also satisfies the
  telescope conjecture. By \cite[Proposition 7.10]{HS}, we then
  have
  \[ L_{K(n-1)}E_n = \holim S/(p^{i_0}, v_1^{i_1}, \dotsb, v_{n-2}^{i_{n-2}})
    \sma v_{n-1}^{-1}E_n, \] where the limit is over type $(n-1)$ generalized
  Moore spectra. We observe that
  \[ (v_{n-1}^{-1}E_n)_*S/(p^{i_0}, v_1^{i_1}, \dotsb, v_{n-2}^{i_{n-2}}) =
    E_*[u_{n-1}^{-1}]/(p^{i_0}, u_1^{i_1}, \dotsb, u_{n-2}^{i_{n-2}}), \] which
  is even periodic with
  \[ (v_{n-1}^{-1}E_n)_0S/(p^{i_0}, v_1^{i_1}, \dotsb, v_{n-2}^{i_{n-2}}) =
    Wk[[u_1, \dotsc, u_{n-1}]][u_{n-1}^{\pm 1}]/(p^{i_0}, u_1^{i_1}, \dotsb,
    u_{n-2}^{i_{n-2}}). \] The transition maps in the diagram are surjective, so
  there is no $\lim\nolimits^1$ and the result is still even periodic. The limit
  on $\pi_0$ is the completion $Wk[[u_1, \dotsc, u_{n-1}]][u_{n-1}^{\pm
    1}]^{\wedge}_{(p, u_1, \dotsc, u_{n-2})}$, as desired.
\end{proof}

\begin{proposition}
  We have
  \[ \LE_0 = Wk((u_{n-1}))^\wedge_p[[u_1, \dotsc, u_{n-2}]]. \]
\end{proposition}

\begin{proof}
  Elements of both rings can be identified as certain possibly infinite formal
  sums
  \[ \sum \{a_Iu_1^{i_1}\dotsm u_{n-2}^{i_{n-2}}u_{n-1}^{i_{n-1}}: a_I \in Wk,
    i_j \in \N\text{ for }1 \le j \le n-2, i_{n-1} \in \Z.\} \] Such a sum is in
  $\LE_0$ if and only if its reduction modulo each power of $(p, u_1, \dotsc,
  u_{n-2})$ is in $k((u_{n-1}))$. In other words, the exponents $i_{n-1}$
  appearing in all nonzero terms with $i_0, \dotsc, i_{n-2}$ less than some
  fixed $i$ are bounded below. On the other hand, such a sum is in
  $Wk((u_{n-1}))^\wedge_p[[u_1, \dotsc, u_{n-2}]]$ if and only if the terms with
  each fixed $i_1, \dotsc, i_{n-2}$ add up to an element of $u_1^{i_{n-1}}\dotsm
  u_{n-2}^{i_{n-2}}Wk((u_{n-1}))^\wedge_p$. That is, the exponents $i_{n-1}$
  appearing in the nonzero terms with fixed $i_1, \dotsc, i_{n-2}$, and with
  $i_0$ less than some fixed $i$, are bounded below. Since there are only
  finitely many choices of $i_1, \dotsc, i_{n-2}$ less than any fixed $i$, the
  two conditions are in fact equivalent.
\end{proof}

Forgetting about the ring structure, the spectrum $\LE$ also has a simple
description.

\begin{proposition}[{\cite[Theorem 3.10]{AMS}}]\label{additive splitting}
  The spectrum $\LE$ splits as a $K(n-1)$-local coproduct of copies of $E_{n-1}$.
\end{proposition}

\subsection{Completed homology}

\begin{definition}
  We write $\CplMod_{\LE_*}$ for the category of graded $\LE_*$-modules which
  are (degreewise) L-complete with respect to the ideal $I_{n-1} = (p, \dotsc,
  u_{n-2})$.
\end{definition}

The spectrum $\LE$ defines completed homology and cohomology theories:
\[ \LE^*X = \pi_*F(X, \LE), \]
\[ \LE^\wedge_*X = \pi_*L_{K(n-1)}(\LE \sma X) = \pi_*L_{K(n-1)}(E \sma X). \]

\begin{proposition}
  The functors $\LE^*$ and $\LE^\wedge_*$ from $\Ho\Top$ to $\Mod_{\LE_*}$
  naturally factor through $\CplMod_{\LE_*}$.
\end{proposition}

\begin{proof}
  The homology of the sphere is complete, and thus L-complete. Since
  $\CplMod_{\LE_*}$ is an abelian category closed under extensions, the same
  follows for any finite complex. Now let $X$ be an arbitrary spectrum and write
  $X$ as a filtered colimit of its finite subcomplexes $X_\alpha$. Then $\LE^*X
  = \lim \LE^* X_\alpha$, which is also L-complete.

  Finally, the completed homology of $X$ is
  \[ \LE_*^\wedge X = \pi_*L_{K(n-1)}(E_n \sma L_{n-1}X) = \pi_*\holim
    (E[u_{n-1}^{-1}]/(p^i_0, \dotsc, u_{n-2}^{i_{n-2}}) \sma L_{n-1}X). \] There
  is a Milnor exact sequence
  \begin{multline*}
    0 \to \lim\nolimits^1 \pi_{k+1}(E[u_{n-1}^{-1}]/(p^i_0, \dotsc, u_{n-2}^{i_{n-2}}) \sma L_{n-1}X) \to \\
    \LE_k^\wedge X \to \lim \pi_k (E[u_{n-1}^{-1}]/(p^i_0, \dotsc,
    u_{n-2}^{i_{n-2}}) \sma L_{n-1}X) \to 0.
  \end{multline*}
  Each term in the limit diagram is torsion to a power of $I_{n-1}$ and thus
  L-complete. Since L-complete modules are closed under extensions, limits and
  $\lim\nolimits^1$, $\LE_k^\wedge X$ is L-complete for each $k$.
\end{proof}

\begin{proposition}
  If $X$ is finite, then $\LE_*^\wedge X = \LE_*X$, which is complete in the
  ordinary sense.
\end{proposition}

\begin{proof}
  If $X$ is finite, then $\LE \sma X$ is in the thick subcategory generated by
  $\LE$. In particular, it is $K(n-1)$-local. It follows that $\LE_* X =
  \LE_*^\wedge X$, a finite L-complete $\LE_*$-module. Since it is finitely
  generated, it is also complete in the ordinary sense \cite{GM}.
\end{proof}

\begin{proposition}
  If $\LE_*X$ is free over $\LE_*$, then $\LE_*^\wedge X$ is its
  $I_{n-1}$-completion.
\end{proposition}

\begin{proof}
  In the Milnor exact sequence for $\LE_*^\wedge X$, we have $\pi_*(\LE/I \sma
  X) \cong \pi_*(\LE \sma X)/I$, where $I$ is an ideal in $\LE_0$. Thus, the
  transition maps in the towers are surjective, the $\lim\nolimits^1$ term
  vanishes, and the $\lim$ term is the ordinary completion of $\LE_*X$.
\end{proof}

\subsection{Augmented deformations}

\begin{definition}
  Write $\CLaur$ for the coefficient ring of completed Laurent series,
  $Wk((u_{n-1}))^\wedge_p$. This is a complete local ring with residue field
  $k((u_{n-1})).$ As we have shown,
  $$LE_0 = Wk((u_{n-1}))^\wedge_p[[u_1, \dotsc, u_{n-2}]].$$
\end{definition}

We will also write $\H^u$ for the base change of the universal deformation
formal group $\G^u$ over $E_0$ to $\LE_0$, and $\H$ for its base change to the
residue field $k((u_{n-1}))$. By the discussion in \Cref{parameters}, if we
started with the Honda formal group law with $p$-series
\[ [p]_\Gamma(x) = x^{p^n}, \] then $\H$ has a coordinate with $p$-series
\[ [p]_\H(x) = u_{n-1} x^{p^{n-1}} +_\H x^{p^n}. \] In particular, its height is
$n-1$.

\begin{definition}
  Let $\H$ be a formal group over $k((u_{n-1}))$. An \textbf{augmented
    deformation} of $\H$ over $(R, \mathfrak{m}) \in \CLN$ is a triple $(\G, i,
  \alpha)$, where:
  \begin{itemize}
  \item $\G$ is a formal group over $R$,
  \item $i:\CLaur \to R$ is a local ring homomorphism (that is, continuous for
    the maximal ideal topology), inducing a map $\overline{i}: k((u_{n-1})) \to
    R/\mathfrak{m}$,
  \item and $\alpha: \H \otimes^{\overline{i}}_{k((u_{n-1}))} R/\mathfrak{m}
    \stackrel{\sim}{\to} \G \otimes_{R} R/\mathfrak{m}$ is an isomorphism of
    formal groups over $R/\mathfrak{m}$.
  \end{itemize}

  An \textbf{isomorphism} of augmented deformations of $\H$ over $R$,
  $\phi:(\G_1, i_1, \alpha_1) \to (\G_2, i_2, \alpha_2)$, is
  \begin{itemize}
  \item the condition that $i_1 = i_2$,
  \item and a map $\phi:\G_1 \to \G_2$ of formal groups over $R$,
  \item such that the square
    \[ \xymatrix{ \Gamma \otimes_{k((u_{n-1}))} R/\mathfrak{m} \ar[d]^1
        \ar[r]^-{\alpha_1} & \G_1 \otimes_R R/\mathfrak{m} \ar[d]_{\phi} \\
        \Gamma \otimes_{k((u_{n-1}))} R/\mathfrak{m} \ar[r]_-{\alpha_2} & \G_2
        \otimes_R R/\mathfrak{m} } \] commutes.
  \end{itemize}

  Let
  \[ \Def_\H^\aug:\CLN \to \Gpd \] be the functor that sends $R$ to the groupoid
  of augmented deformations of $\H$ over $R$ and their isomorphisms.
\end{definition}

\begin{thm}\label{LE representability}
  The functor $\Def_\H^\aug$ is represented by $\LE_0$. That is, for $R$ a
  complete local ring, the groupoid $\Def_\H^\aug(R)$ is naturally equivalent to
  the discrete groupoid $\Maps_{cts}(\LE_0,R)$ of continuous maps with respect
  to the maximal ideal topology on both rings.
\end{thm}

\begin{proof}
  This is more or less an immediate consequence of the Lubin-Tate theorem, in
  the form of \Cref{Lubin-Tate}. Consider the following deformation in
  $\Def^\aug_\H(\LE_0)$:
  \[
    (\H^u, \,\, 1:\CLaur \to \CLaur, \,\, \alpha^u:\H \stackrel{\sim}{\to} \H^u
    \otimes k((u_{n-1}))),
  \]
  where $\alpha^u$ is the canonical isomorphism given by the definition of $\H$.
  We will show that this deformation is universal.

  Given a local ring homomorphism $f: LE_0 \to R$, we obtain an augmented
  deformation
  \[
    (\H^u \otimes_{LE_0}^f R, \,\, f|_\CLaur, \,\, \alpha^u \otimes^f_{LE_0} R)
    \in \Def^\aug_\H(R)
  \]
  by base change. Note that the map $f|_\CLaur$ is local because $f$ is.

  On the other hand, suppose given $(\G, i, \alpha) \in \Def^\aug_\H(R)$. We
  must exhibit a unique continuous map $f:\LE_0 \to R$ and a unique isomorphism
  between $(\G, i, \alpha)$ and an augmented deformation of the above form.
  Since $i$ is local, we may regard $R$ as a $\CLaur$-algebra of the form given
  in \Cref{Lubin-Tate}, and $(\G, \overline{i}, \alpha)$ as an object of
  $\Def_\H^\CLaur(R)$. \Cref{Lubin-Tate} now implies that $(\G, i, \alpha)$ is
  uniquely isomorphic, as an object of $\Def_\H^\CLaur(R)$, to a base change of
  the universal deformation along a unique local ring map
  \[ \CLaur[[u_1, \dotsc, u_{n-2}]] \to R \] compatible with $i$. Equivalently,
  it is uniquely isomorphic, as an object of $\Def_\H^\aug(R)$, to a base change
  of $(\H^u, 1, \alpha^u)$ along a unique local ring map
  \[ \CLaur[[u_1, \dotsc, u_{n-2}]] \to R. \]
\end{proof}

\section{The \texorpdfstring{$E$}{E}-theory of
  \texorpdfstring{$E$}{E}-theory}\label{FE}

In this section, we let $k$ be a finite field containing $\F_{p^n}$ and
$\F_{p^{n-1}}$. We let $E$ be the $E$-theory associated to a height $n$ formal
group $\Gamma$ over $k$, and $F$ the $E$-theory associated to the height $n-1$
Honda formal group $\Gamma_{n-1}$ over $k$. We have just described the homotopy
groups of $\LE$; as $\LE$ is a $K(n-1)$-local spectrum, it is natural to study
it in terms of its completed $F$-homology, which we do in this section. We will
also obtain a description of the coalgebra of cooperations, $\LE^\wedge_* \LE$.

\subsection{\texorpdfstring{The completed $E_{n-1}$-homology of $E_n$}{The
    completed En-1-homology of En}}

\begin{proposition}
  $F_*^\wedge E = F_*^\wedge \LE$.
\end{proposition}

\begin{proof}
  The map $E \to \LE$ is a $K(n-1)$-local equivalence, so it remains so after
  smashing with $F$.
\end{proof}

\begin{proposition}\label{FE flat}
  $F_*^\wedge E$ is even periodic and flat over $F_*$.
\end{proposition}

\begin{proof}
  As in the proof of \Cref{E cooperations periodic}, $F$ and $\LE$ are both even
  periodic and Landweber exact, so
  \[ F_* \LE = F_* \otimes_{BP_*} BP_*BP \otimes_{BP_*} \LE_*, \] which is even
  periodic and flat over $F_*$, since $F \sma E$ is $L_{n-1}$-local. The
  $K(n-1)$-localization satisfies
  \[ F_*^\wedge \LE = \pi_*\holim_I (F \sma \LE/I), \] where $I$ ranges over a
  cofinal set of ideals of the form $(p^{i_0}, \dotsc, v_{n-1}^{i_{n-1}})$. The
  objects in the limit diagram are
  \[ F_* \otimes_{BP_*} BP_*BP \otimes_{BP_*} \LE_*/I, \] and the maps in the
  diagram are surjective. Therefore, $F_*^\wedge E$ is also even periodic.

  It remains to show that $F_*^\wedge E$ is $F_*$-flat. The above implies that
  \[ F_*^\wedge \LE = (F_* \LE)^\wedge_{I_{n-1}}, \] (the degreewise
  completion), and therefore that
  \[ F_*^\wedge \LE/I = (F_*\LE)/I = F_*/I \otimes_{BP_*} BP_*BP \otimes_{BP_*}
    \LE_*/I, \] for any power $I$ of $I_{n-1}$. Since each of these is a flat
  $F_*/I$-module and the maps in the limit diagram computing the completion are
  surjective, \cite[Tag 0912]{Stacks} implies that $F_*^\wedge \LE$ is a flat
  $(F_*)^\wedge I = F_*$-module.
\end{proof}

We now describe the functor represented by $F_0^\wedge E$.

\begin{lemma}
  For any complete local ring $R$, pre-composition with the completion map $F_0E
  \to F_0^\wedge E$ induces an isomorphism
  \[ \Hom_{cts}(F_0^\wedge E, R) \cong \Hom_{F_0,cts}(F_0E, R). \] Moreover, for
  any continuous map $F_0^\wedge E \to F_0$, the composition
  \[ \LE_0 \stackrel{\eta_R}{\to} F_0E \to R \] is also continuous.
\end{lemma}

\begin{proof}
  This is just as in \Cref{E cooperations continuity}. Since $I_{n-1}$ is an
  invariant ideal in $BP_*$, any complete local ring $R$ with a map $f:F_0E \to
  R$ such that the restriction to $F_0$ is continuous must be complete with
  respect to $I_{n-1}\cdot F_0E$, so that $f$ factors uniquely through the
  completion $F^\wedge_0 E$. Moreover, $f$ also sends $\eta_R(I_{n-1}\cdot E_0)$
  into $I_{n-1} R$, so that the map $E_0 \to R$ coming from the right unit is
  also continuous.
\end{proof}

\begin{thm}\label{FE theorem}
  Let $R$ be a complete local $F_0$-algebra. There is a natural isomorphism
  between continuous $F_0$-algebra maps $F_0^\wedge E \to R$ and pairs
  $(j,\gamma)$, where $j:\CLaur \to R$ is a $p$-adically continuous map and
  $\gamma$ is an isomorphism of formal groups over $R/\mathfrak{m}$,
  $\gamma:\Gamma_{n-1} \otimes_k^{i} R/\mathfrak{m} \stackrel{\sim}{\to} \H
  \otimes_{k((u_{n-1}))}^{\ol{j}} R/\mathfrak{m}$.
\end{thm}

\begin{proof}
  As before, we have
  \[ F_0\LE = \pi_0(F_* \otimes_{BP_*} BP_*BP \otimes_{BP_*} \LE_*) = F_0
    \otimes_{BPP_0} BPP_0BPP \otimes_{BPP_0} \LE_0. \] An $F_0$-algebra map
  $F_0\LE \to R$ is equivalent to a map $\LE_0 \to R$ and an isomorphism over
  $R$ between the base changes of the formal groups of $F_0$ and $\LE_0$.

  If $R$ is complete local, then the previous lemma tells us that
  $\Hom_{F_0,cts}(F_0E,R) = \Hom_{F_0,cts}(F_0^\wedge E, R)$, and that the map
  $\LE_0 \to R$ is continuous. Hence, the map $\LE_0 \to R$ represents an object
  of $\Def_\H^\aug(R)$. Likewise, the structure map $F_0 \to R$ represents an
  object of $\Def_{\Gamma_{n-1}}(R)$, say $(\G, i, \alpha)$. Thus, we have a
  pullback of groupoids:
  \[ \xymatrix{ \Hom_{F_0,cts}(F_0^\wedge E, R) \ar[r] \ar[d] \pullbackcorner &
      \Def_\H^\aug(R) \ar[d] \\ \{*\} \ar[r]_{\G} & \Mfg(R). } \] In other
  words, a map $f:F_0E \to R$ corresponds to the data:
  \begin{gather*} (\G', j:\Lambda \to R, \beta:\H
    \otimes_{k((u_{n-1}))}^{\ol{j}} R/\mathfrak{m} \stackrel{\sim}{\to} \G'
    \otimes_R R/\mathfrak{m}) \in
    \Def_\H^\aug(R); \\
    \phi:\G \stackrel{\sim}{\to}{\G'}.
  \end{gather*}
  There is a unique isomorphism in the pullback groupoid which restricts to
  $\phi^{-1}:\G' \to \G$ on formal groups. Composing with this isomorphism, one
  gets a unique object in the pullback groupoid isomorphic to $f$ whose
  underlying formal group is $\G$ and whose underlying isomorphism of formal
  groups is the identity. It follows that the groupoid is locally contractible,
  as expected. The rest of the data is given by $j$ and $\beta$, or equivalently
  by $j$ and
  \[ \gamma = \beta^{-1}\alpha:\Gamma_{n-1} \otimes_k^{i} R/\mathfrak{m}
    \stackrel{\sim}{\to} \H \otimes_{k((u_{n-1}))}^{\ol{j}} R/\mathfrak{m}. \]
\end{proof}

Now let's see what this functorial description of $F^\wedge_*E$ can tell us
about it algebraically.

\begin{proposition}\label{FE finite product}
  $F_0 E/I_{n-1}$ is of the form $\Hom(\Gal(k/\F_p),L)$, where $L$ is a field.
  Therefore, $F_0^\wedge E$ is a finite product of complete local rings.
\end{proposition}

\begin{proof}
  Armed with \Cref{FE theorem}, this is essentially a reinterpretation of a
  result of Torii, \cite[Theorem 2.7]{Torii}, which in turn reinterprets a
  result from \cite{Gross}. For $R$ a complete local $k = F_0/I_{n-1}$-algebra,
  \[ \Hom_{k}(F_0E/I_{n-1}, R) = \{(\ol{j}:k((u_{n-1})) \to R,
    \gamma:\Gamma_{n-1} \otimes_k^{i} R/\mathfrak{m} \stackrel{\sim}{\to} \H
    \otimes_{k((u_{n-1}))}^{\ol{j}} R/\mathfrak{m})\}. \] The \'etaleness of
  isomorphisms, \Cref{etaleness of isomorphisms}, says that we can equivalently
  define $\gamma$ as an isomorphism between $\Gamma_{n-1}$ and $\H$ over $R$.
  There is a smallest extension $L$ of $k((u_{n-1}))$ over which $\Gamma_{n-1}$
  and $\H$ become isomorphic, given by adjoining the coefficients of an
  isomorphism between any choice of formal group laws for $\Gamma_{n-1}$ and
  $\H$. Torii proves that $L/k((u_{n-1}))$ is Galois with Galois group
  $\Aut_k(\Gamma_{n-1})$. On the other hand, having chosen $\ol{j}:k((u_{n-1}))
  \to R$ such that an isomorphism $\Gamma_{n-1} \otimes R \to \H \otimes R$
  exists, the set of such isomorphisms is clearly a torsor for this group. It
  follows that a $k$-algebra map $F_0E/I_{n-1} \to R$ is precisely a map $L \to
  R$, which is not necessarily a $k$-algebra map. In other words, $F_0E/I_{n-1}
  = \Hom(\Gal(k/\F_p),L)$.

  Since $k$ is finite, this is a finite product of fields. The corresponding
  splitting for $F_0^\wedge E$ itself follows from Hensel's lemma.
\end{proof}

\begin{remark}
  It's possible to be slightly more explicit, using the formula
  \[
    F_*^\wedge E/I_{n-1} = F_* \otimes_{BP_*} BP_*BP \otimes_{BP_*}
    LE_*/I_{n-1}.
  \]
  For $x \in BP_*$, write $x$ for $\eta_L(x)$ and $\overline{x}$ for
  $\eta_R(x)$, elements of $BP_*BP$. To avoid confusion, let $u_i$ and $u$ be
  the Lubin-Tate generators of $E_*$, and write $u_{i}^F$ and $u^F$ for the
  corresponding generators in $F_*$. The map $BP_* \to F_*$ sends
  \begin{align*}
    v_i &\mapsto (u^F)^{p^i-1}u_i^F, \quad i \le n-2, \\
    v_{n-1} &\mapsto (u^F)^{p^{n-1} - 1}, \\
    v_i &\mapsto 0, \quad i \ge n.
  \end{align*}
  Likewise, the map $BP_* \to \LE_*$ sends
  \begin{align*}
    \ol{v_i} &\mapsto u^{p^i-1}u_i, \quad i \le n-1, \\
    \ol{v_n} &\mapsto u^{p^{n} - 1}, \\
    \ol{v_i} &\mapsto 0, \quad i \ge {n+1}
  \end{align*}

  We can thus write
  \[ F_*^\wedge E = \left(\frac{Wk[[u_1^F, \dotsc, u_{n-2}^F]][(u^F)^{\pm
          1}][t_1, t_2, \dotsc] \otimes_{\F_p} Wk((u_{n-1}))[[u_1, \dotsc,
        u_{n-2}]][u^{\pm 1}]}{(\ol{v_1} - u^{p^i-1}u_i, \dotsc, \ol{v_n} -
        u^{p^{n} - 1}, \ol{v_{n+1}}, \dotsc)}\right)^\wedge_{I_{n-1}}. \] In
  degree zero, let $s_i = t_iu^{1 - p^i}$ and $w = u^F/u$. Note that the ideal
  $I_{n-1}$ contains $p$ and all $u_i$ and $u_i^F$ (and thus all $v_i$ and
  $\ol{v_i}$) for $1 \le i \le n-2$. Therefore,
  \[ F_0^\wedge E/I_{n-1} = (k \otimes_{\F_p} k)((u_{n-1}))[s_1, s_2, \dotsc,
    w^{\pm 1}]/(u^{1-p^{n-1}}\ol{v_{n-1}} - u_{n-1}, u^{1-p^n}\ol{v_{n}} - 1,
    \ol{v_{n+1}}, \dotsc). \] Now, by \cite[4.3.1]{Rav},
  \begin{equation}\label{BP right unit}
    \ol{v_{n-1+i}} \equiv v_{n-1+i} + v_{n-1}t_{i}^{p^{n-1}} - v_{n-1}^{p^{i}}t_{i}\pmod{(p, v_1, \dotsc, v_{n-2}, t_1, \dotsc, t_{i-1})}.
  \end{equation}
  Scaling to degree zero by multiplying by appropriate powers of $u$, we get
  relations in $F_0^\wedge E/I_{n-1}$:
  \begin{align}\label{si relations}
    \begin{split}
      u_{n-1} &= w^{p^{n-1}-1}, \\
      1 &= w^{p^{n-1}-1}s_1^{p^{n-1}} - w^{p(p^{n-1}-1)}s_1, \\
      0 &\equiv w^{p^{n-1}-1}s_i^{p^{n-1}} - w^{p^i(p^{n-1}-1)}s_i + f_i(s_1,
          \dotsc, s_{i-1})\text{ for }i \ge 2.
    \end{split}
  \end{align}
  The first relation gives an embedding of $k((w))$ into $F_0^\wedge E/I_{n-1}$,
  as a tamely ramified extension of $k((u_{n-1}))$ (recall that $k$ contains
  $\F_{p^{n-1}}$). The remaining relations inductively define $s_i$ as solutions
  to higher Artin-Schreier equations over the ring generated by $k((w))$ and
  $s_1, \dotsc, s_{i-1}$. In particular, $F_0^\wedge E/I_{n-1}$ is ind-\'etale
  over $k((u_{n-1}))$.

  What is not clear from this calculation, but follows from Torii's theorem, is
  that the higher Artin-Schreier equations in \eqref{si relations} are
  irreducible, thus defining a field extension $L/k((u_{n-1}))$.
\end{remark}

\begin{remark}
  One should compare this result with the computation of $E^\wedge_0 E$. In both
  cases, the object calculated is a flat extension of a Lubin-Tate ring, it
  represents an isomorphism of deformations of formal groups, and its reduction
  mod $I_{n-1}$ carries the action of a Morava stabilizer group. In the case of
  $E^\wedge_0 E$, this action splits, and in fact $E^\wedge_0 E$ is a profinite
  group algebra for the Morava stabilizer group. In the case of $F_0^\wedge E$,
  the action is in a certain sense as complicated as possible, so that the only
  part of the action that splits is the Galois group. By Torii's theorem, the
  other part $\Aut_k(\Gamma_{n-1})$ of the Morava stabilizer group instead
  becomes the Galois group of a residue field extension $L/k((u_{n-1}))$. In
  particular, $\H$ only becomes isomorphic to the Honda formal group
  $\Gamma_{n-1}$ after base change to $L$. This is one way of saying that the
  formal group of $\LE_0$ is as complicated a height $n-1$ formal group as
  possible.
\end{remark}

The formal \'etaleness of $E_0^\wedge E_0$ over $E_0$ is measured by the
vanishing of its completed cotangent complex. (To be precise, this is the
complex representing derived functors of derivations of $E_0$-algebras which are
L-complete with respect to $I_n$.) $F_0^\wedge E$ is not formally \'etale over
$F_0$, but it is formally smooth, as shown by the following calculation. (The
cotangent complex in question is likewise defined using the category of
$F_0$-algebras which are L-complete with respect to $I_{n-1}$.)

\begin{proposition}\label{formally smooth}
  The completed cotangent complex $\L_{F_0^\wedge E/F_0}$ is concentrated in
  degree zero and free of infinite rank (equal to the transcendence degree of
  $k((u_{n-1}))/k$).
\end{proposition}

\begin{proof}
  In \Cref{FE theorem}, we showed that $F_0^\wedge E$ represents data of the
  following form on complete $F_0$-algebras:
  \[
    (j:\CLaur \to R, \gamma\text{ an isomorphism of formal groups over }
    R/\mathfrak{m}).
  \]
  But the $\gamma$ part of the data is insensitive to infinitesimal thickenings.
  Thus, if $R$ is a complete local ring and $I$ a square-zero ideal,
  constructing a lift in the diagram
  \[ \xymatrix{ \LE_0 \ar[r] \ar[d] & R \ar[d] \\ \LE_0^\wedge \LE \ar[r]
      \ar@{-->}[ur] & R/I } \] reduces to constructing a lift in the diagram
  \[ \xymatrix{ Wk \ar[r] \ar[d] & R \ar[d] \\ \CLaur = Wk((u_{n-1}))^\wedge_p
      \ar[r] \ar@{-->}[ur] & R/I. } \] In other words, the complete cotangent
  complex $\L_{\LE_0^\wedge \LE/\LE_0}$ is a base change of $\L_{\CLaur/Wk}$. As
  a consequence of \Cref{existence of cohen rings}, the map $Wk \to \CLaur$ is
  formally smooth if $k$ is perfect, meaning that its completed cotangent
  complex is $\CLaur$-free and concentrated in degree zero, and
  \[ \rank_\CLaur \L_{\CLaur/Wk} = \rank_{k((u_{n-1}))} \L_{k((u_{n-1}))/k} \]
  which is equal to the transcendence degree of $k((u_{n-1}))/k$. (This is
  always infinite, but there is an easy cardinality argument if $k$ is finite:
  $k((u_{n-1}))$ is uncountable, while any field of finite transcendence degree
  over $k$ is countable.)
\end{proof}

Finally, we make a comment about the $K(n-1)$-local Adams-Novikov spectral
sequence for $\LE$. Let $\G_{n-1} = \Aut(k, \Gamma_{n-1})$ be the height $n-1$
Morava stabilizer group. By \Cref{E cooperations} and the fact that $F_*^\wedge
E$ is honestly complete, we can regard its $F_*^\wedge F$-comodule structure as
a continuous $\G_{n-1}$-module structure.

\begin{lemma}\label{action extended}
  The $\G_{n-1}$-action on $F_*^\wedge E$ is extended. In other words, there is
  a $\G_{n-1}$-module isomorphism:
  \[
    F_*^\wedge E \cong \Cts(\G_{n-1}, \LE_*).
  \]
\end{lemma}

\begin{proof}
  This is true for $F$, and $\LE$ is a $K(n-1)$-local sum of copies of $F$, by
  \Cref{additive splitting}
\end{proof}

\begin{corollary}\label{ANSS collapses}
  The $K(n-1)$-local Adams-Novikov spectral sequence
  \[
    E_2 = H^p_{cts}(\G_{n-1}, F_*^\wedge E) \Rightarrow \pi_*\LE
  \]
  is concentrated on the 0-line and collapses at $E_2$:
  \[
    \pi_*\LE = (F_*^\wedge E)^{\G_{n-1}}
  \]
\end{corollary}

\subsection{Cooperations for \texorpdfstring{$LE$}{LE}}

The same arguments apply to the completed cooperations algebra $\LE_*^\wedge
\LE$. Let us merely state the corresponding results without proof:

\begin{proposition}\label{LE cooperations}
  \hfill
  \begin{enumerate}[(a)]
  \item $LE_*^\wedge \LE$ is even periodic and flat over $\LE_*$.
  \item Let $R$ be a complete local $\LE_0$-algebra, with $\LE_0 \to R$
    classifying $(\G, i, \alpha) \in \Def_\H^{\aug}(R)$. There is a natural
    isomorphism between continuous $\LE_0$-algebra maps $\LE_0^\wedge \LE \to R$
    and pairs $(j,\gamma)$, where $j:\CLaur \to R$ is a $p$-adically continuous
    map and $\gamma$ is an isomorphism of formal groups over $R/\mathfrak{m}$,
    $\gamma:\H \otimes_{k((u_{n-1}))}^{\ol{i}} R/\mathfrak{m}
    \stackrel{\sim}{\to} \H \otimes_{k((u_{n-1}))}^{\ol{j}} R/\mathfrak{m}$.
  \item The completed cotangent complex $\L_{LE_0^\wedge LE/LE_0}$ is
    concentrated in degree zero and free of infinite rank (equal to the
    transcendence degree of $k((u_{n-1}))/k$.)
  \end{enumerate}

\end{proposition}

\section{Power operations and \texorpdfstring{$\Einf$}{Einfty}
  structures}\label{power ops}

\subsection{Background on \texorpdfstring{$\theta$}{theta}-algebras}

This section specializes to the case $n = 2$ and $n-1 = 1$. We write $K = E_1$
and $E = E_2$, both over a finite field $k$ containing $\F_{p^2}$. We write $x$
for $u_1$, so that $LE_0 \cong Wk((x))^\wedge_p$.

\begin{definition}
  A \textbf{$\theta$-algebra} is a $\Z_p$-algebra $R$ equipped with an operation
  $\theta:R \to R$ such that
  \begin{align*}
    \theta(x+y) &= \theta(x) + \theta(y) - \sum_{i=1}^{p-1} \frac{1}{p}\binom{p}{i} x^i y^{p-i}, \\
    \theta(xy) &= \theta(x)y^p + x^p\theta(y) + p\theta(x)\theta(y), \\
    \theta(0) &= \theta(1) = 0.
  \end{align*}
\end{definition}

These properties imply that the operation $\psi^p:R \to R$ defined by
\[ \psi^p(x) = x^p + p\theta(x) \] is a ring homomorphism; in other words, any
$\theta$-algebra $R$ is equipped with a lift of the Frobenius on $R/p$. If $R$
is $p$-torsion-free, then the Frobenius lift $\psi^p$ uniquely determines
$\theta$, but in general, $\theta$ is strictly more data than $\psi^p$.

\begin{definition}
  A \textbf{$\theta$-comodule algebra} is a $\theta$-algebra $R$ together with a
  continuous action of $\Z_p^\times$ on $R$ that commutes with $\theta$. We
  write $\ComodAlg_\theta$ for the category of $\theta$-comodule algebras.
\end{definition}

In $K(1)$-local homotopy theory, $\theta$-algebras arise from power operations
on $\Einf$ ring spectra. More precisely:
\begin{itemize}
\item If $X$ is a $K(1)$-local $\Einf$ ring spectrum, then $K^\wedge_0X$ is
  naturally a $\theta$-comodule algebra \cite[Theorem 2.2.4]{GH2}.
\item If $X$ is a $K(1)$-local $\Einf$ ring spectrum, then $\pi_0X$ is naturally
  a $\theta$-algebra \cite{HopkinsK1local}.
\end{itemize}

In the remainder of this section, we consider the case $X = \LE$. In section
6.2, we show that there are nonisomorphic $\theta$-algebra structures on
$\LE_0$. In section 6.3, we extend these structures to $\theta$-comodule algebra
structures on $K_0^\wedge E$. In section 6.4, we use an obstruction theory
argument to show that any $\theta$-algebra structure on $\LE_0$ is realized by
some $\Einf$ structure on $\LE$.

\subsection{Constructing non-isomorphic \texorpdfstring{$\theta$}{theta}-algebra
  structures}

Let $X$ and $X'$ be two $\Einf$-algebras abstractly equivalent to $LE$. An
equivalence $X \to X'$ induces a $p$-adically continuous isomorphism of
$\theta$-algebras $\pi_0X \to \pi_0X'$. Thus, the question of identifying
$\Einf$ structures on $LE$ up to equivalence is related to the question of
identifying $\theta$-algebra structures on $LE_0 = Wk((x))^\wedge_p$ up to
isomorphism.

I thank Dominik Absmeier for catching a mistake in the original version of this
argument.

\begin{lemma}\label{theta mod p}
  Given any $\theta$-algebra structure on $Wk((x))^\wedge_p$, the operation
  $\theta$ descends to a map $\theta:W_2k((x)) \to k((x))$.
\end{lemma}

\begin{proof}
  We have to show that $\theta(f)$ mod $p$ only depends on the class of $f$ mod
  $p^2$. First, because $\psi$ is a ring homomorphism, we must have
  \[
    \psi(n) = n\quad\text{for }n \in \Z;
  \]
  this forces
  \[
    \theta(n) = \frac{n - n^p}{p}\quad\text{for }n \in \Z.
  \]
  In particular, if $v_p(n) = i \ge 1$, then $v_p(\theta(n)) = i-1$.

  Now calculate
  \begin{align*}
    \theta(f + p^2g) &= \theta(f) + \theta(p^2g) - \sum_{i=1}^{p-1}
                       \frac{1}{p}\binom{p}{i} f^i p^{2(p-i)}g^{2(p-i)} \\
                     &= \theta(f) + p^{2p}\theta(g) + \theta(p^2)g^{p} + p\theta(p^2)g^p - O(p^2) \\
                     &= \theta(f) + O(p)
  \end{align*}
  because $\theta(p^2)$ is divisible by $p$.
\end{proof}

We will show that there are two $\theta$-algebra structures, $\theta_0$ and
$\theta$, such that there is no $p$-adically continuous automorphism of
$Wk((x))^\wedge_p$ making the diagram
\[ \xymatrix{ Wk((x))^\wedge_p \ar[r]^{\theta_0} \ar[d]_f & Wk((x))^\wedge_p
    \ar[d]_f \\ Wk((x))^\wedge_p \ar[r]_{\theta} & Wk((x))^\wedge_p } \]
commute.

To do this, it helps to introduce an alternate notion of continuity, besides the
$p$-adic one we have been using to this point. We first give some motivation,
recalling \Cref{three continuities} from the introduction. The ring
$Wk((x))^\wedge_p = Wk((u_{n-1}))^\wedge_p$ arose in the first place via
inverting an element in the maximal ideal of a complete local ring and then
completing with respect to a smaller ideal -- these algebraic operations being
related to the topological one of $K(n-1)$-localizing a $K(n)$-local object.
Treating $Wk((x))^\wedge_p$ as a $p$-adically complete local ring, equipped with
its $p$-adic topology, means forgetting about the topology that was previously
on the ring $Wk[[x]]$. This gives us a lot of freedom in specifying continuous
maps out of $Wk((x))^\wedge_p$ -- just as there are many more discontinuous than
continuous maps out of the ring $k[[x]]$ -- but also makes them harder to
specify -- just as a continuous ring map out of $k[[x]]$ is determined by the
image of $x$, but a discontinuous one is not.

Instead, one could keep track of both the $p$-adic topology on
$Wk((x))^\wedge_p$ and the maximal ideal topology on $Wk[[x]]$. This is
formalized in the notion of \textbf{pipe structure} due to \cite{MGPS}, drawing
on work on higher local fields in number theory \cite{Kato, Morrow}. Rather than
give the full definition, we will describe what it means for a map from
$Wk((x))^\wedge_p$ to itself to be continuous for the pipe structure, as such
maps are all we care about here.

\begin{definition}
  Let $\phi:Wk((x))^\wedge_p \to Wk((x))^\wedge_p$ be a map of sets (which is
  not necessarily a ring homomorphism). We say that $\phi$ is
  \textbf{pipe-continuous} if it is the limit of an inverse system of maps
  \[ \phi_i: W_{n_i}k((x)) \to W_{m_i}k((x)), \] where each $W_nk((x))^\wedge_p$
  has the $x$-adic topology.
\end{definition}

\begin{remark}
  The $x$-adic topology on $W_nk((x))$ is the same as the topology induced by
  the maximal ideal topology on $W_nk[[x]]$. Indeed,
  \[ (p, x)^{n+r}W_nk[[x]] \subseteq x^rW_nk[[x]] \subseteq (p,
    x)^rW_nk[[x]]. \]
\end{remark}

To demonstrate the usefulness of this idea, we note that pipe-continuous ring
homomorphisms are very easy to describe.

\begin{proposition}
  A pipe-continuous ring homomorphism $\phi:Wk((x))^\wedge_p \to
  Wk((x))^\wedge_p$ is completely determined by the image of $x$, which can be
  any completed Laurent series that reduces mod $p$ to a nonzero element of
  $xk[[x]]$.
\end{proposition}

\begin{proof}
  First, it is easy to see that a pipe-continuous ring homomorphism $\phi$ can
  be written as the limit of an inverse system
  \[ \phi_i: W_{n_i}k((x)) \to W_{n_i}k((x)), \] and must in particular be
  $p$-adically continuous. Now, any continuous map $W_nk((x)) \to W_nk((x))$ is
  determined by the image of $x$: namely, $x$ can go to any invertible,
  topologically nilpotent element. Let $y \in W_nk((x))$ be invertible and
  topologically nilpotent. Then $y$ is not divisible by $p$, so writing
  \[ y = \sum_{i \ge -N} a_i x^i, \] there is a least $i$, say $i = d$, such
  that $a_d$ is invertible in $W_nk$. If $d \le 0$, then for any $r$, $y^{p^r} =
  \ol{a_d}^{p^r} x^{ip^r} + \dotsb$ mod $p$, so $y$ is not topologically
  nilpotent.

  Conversely, suppose that $y = \sum_{i \ge -n} a_i x^i$, and that $y$ reduces
  mod $p$ to $\ol{a_d}x^d + \dotsb$ with $\ol{a_d} \ne 0$ and $d > 0$. Then mod
  $p$,
  \[ a_d^{-1}x^{-d}y \equiv 1 + \dotsb \in k((x))^\times, \] so that
  $a_d^{-1}x^{-d}y$ is invertible by Hensel's lemma, and so $y$ itself is
  invertible. Now let $y_+$ be the sum of the terms of $y$ of degree $\ge d$, so
  that we can write $y = y_+ + py_-$. Then
  \[ y^{p^r} = y_+^{p^r} + \sum_{i=1}^{p^r} p^i\binom{p^r}{i} y_-^i y_+^{p^r -
      i}. \] By Kummer's theorem on valuations of binomial coefficients,
  \[ v_p\left( p^i\binom{p^r}{i}\right) = i + r - v_p(i) \ge r. \] Thus, for $r$
  sufficiently large, $y^{p^r} = y_+^{p^r}$ in $W_nk((x))$. Thus, $y$ is
  topologically nilpotent.
\end{proof}

\begin{proposition}\label{continuous automorphism}
  Every field automorphism of $k((x))$ is continuous.
\end{proposition}

\begin{proof}
  Suppose that $k = \F_q$. Let $S$ be the set of power series $1 + a_1x +
  \dotsb$. Then $S$ is multiplicatively closed, and any $f \in S$ has a
  $(q-1)$th root $g \in S$, by the binomial theorem. On the other hand, an $f
  \not\in S$ either has $v_{x}(f) = 0$ but constant term not equal to 1, in
  which case it does not have a $(q-1)$th root at all, or $v_{x}(f) \ne 0$, in
  which case it has at most a $(q-1)^m$th root for some maximal $m$.

  It follows that $S$ is exactly the set of elements of $k((x))$ which have a
  $(q-1)^m$th root for all $m$. Thus, any automorphism of $k((x))$ preserves
  $S$. Subtracting 1, we see that any automorphism of $k((x))$ preserves the set
  $xk[[x]]$, and thus that it preserves $x^rk[[x]]$ for every $r$. Thus, any
  automorphism is continuous.
\end{proof}

\begin{corollary}\label{theta continuity}
  If $\psi^p$ is any pipe-continuous Frobenius lift of $Wk((x))^\wedge_p$, then
  the associated $\theta$ induces an ($x$-adically) continuous map
  \[
    \theta: W_2k((x)) \to k((x)).
  \]
\end{corollary}

\begin{proof}
  We can recover this reduction of $\theta$ from $\psi^p$ mod $p^2$, using the
  formula $\psi^p(f) = f^p + p\theta(f)$ inside $W_2k((x))$. In
  other words, $\theta(f) = \frac{1}{p}(\psi^p(f) - p\theta(f))$, where division
  by $p$ is the obvious isomorphism from the $p$-torsion of $W_2k((x))$ to
  $k((x))$. Since $\psi^p(f)$ and $f \mapsto f^p$ are both continuous on
  $W_2k((x))$, it follows that $\theta$ is also continuous.
\end{proof}

\begin{proposition}\label{nonisomorphic theta}
  There exist non-isomorphic $\theta$-algebra structures on $Wk((x))^\wedge_p$.
\end{proposition}

\begin{proof}
  Let $\psi^p_0$ and $\psi^p$ be the unique pipe-continuous endomorphisms
  satisfying $\psi^p_0(x) = x^p$ and $\psi^p(x) = x^p + p$. Thus, $\theta_0(x) =
  0$ and $\theta(x) = 1$. Suppose that $f$ is an automorphism of
  $Wk((x))^\wedge_p$ such that $f\theta_0 = \theta f$. By \Cref{theta mod p},
  there is a commutative diagram
  \[
    \xymatrix{ W_2k((x)) \ar[r]^{\theta_0} \ar[d]_f & k((x)) \ar[d]_{\ol{f}} \\
      W_2k((x)) \ar[r]_{\theta} & k((x)) }
  \]
  where we have written $f$ for the reduction of $f$ mod $p^2$, and $\ol{f}$ for
  the reduction of $f$ mod $p$. In particular, note that
  \[
    \theta(f(x)) = \ol{f}(\theta_0(x)) = 0.
  \]
  We will show that this is not possible for any automorphism $f$.

  By \Cref{continuous automorphism}, $\ol{f}$ is a continuous automorphism of
  $k((x))$. This means that, mod $p^2$, $f$ must be of the form
  \[
    f(x) \equiv ax + g(x) + ph(x^{-1}) \pmod{p^2},
  \]
  where $a \in W_2k^\times$, $g \in x^2W_2k[[x]]$, and $h$ is a polynomial.

  Now let us calculate $\theta(f)$ mod $p$.
  \begin{align*}
    \theta(f(x)) &\equiv \theta(ax) + \theta(g + ph) - \sum \frac{1}{p}\binom{p}{i} x^i(g + ph)^{p-i} \\
                 &\equiv a^p\theta(x) + \theta(a)x^p + \theta(g + ph) - \sum \frac{1}{p}\binom{p}{i}x^ig^{p-i} \\
                 &\equiv a^p + \theta(a)x^p + \theta(g) + \theta(ph) - \sum \frac{1}{p}\binom{p}{i}x^ig^{p-i} \\
                 &\equiv a^p + \theta(a)x^p + \theta(g) + \theta(p)h(x^{-1})^p - \sum \frac{1}{p}\binom{p}{i}x^ig^{p-i}.
  \end{align*}
  Since the only terms of negative degree in $x$ come from
  $\theta(p)h(x^{-1})^p$, we must have $h \equiv 0$ mod $p$, so that the term
  $\theta(p)h(x^{-1})^p$ disappears.

  As for the other terms, write $v_x$ for the $x$-adic valuation on
  $k((x))$. The term $a^p$ has $v_x = 0$ -- thus, for $\theta(f)$ to equal zero,
  it is necessary for some other term to also have $x$-adic valuation 0.
  However,
  \begin{align*}
    v_x(a^p) &= 0, \\
    v_x(x\theta(a)) &= 1, \\
    v_x(x^i g^{p-i}) &= i + (p-i)v_x(g) \ge p+1 \ge 3. \\
  \end{align*}
  As for $g$, still working mod $p$,
  \[
    \theta(x^n) = x^p\theta(x^{n-1}) + x^{p(n-1)}\theta(x) = nx^{p(n-1)},
  \]
  by induction on $n$. In particular,
  \[
    v_x(\theta(x^n)) \ge 2 \text{ when } n \ge 2.
  \]
  By the multiplication formula, the same is true for $\theta(bx^n)$, when $b
  \in W_2k$. Thus, a monomial $g$ with $v(g) \ge 2$ has $v(\theta(g)) \ge 2$ as
  well. By induction and using the theta sum formula, the same is true for
  polynomials. By continuity of $\theta$ (\Cref{theta continuity}), the same is
  true for power series.

  It follows that $v_x(\theta(g)) \ge 2$. Thus, it's impossible for
  $\theta(f(x))$ to equal 0, which is what we needed to prove.
\end{proof}

\subsection{The \texorpdfstring{$K$}{K}-theory of \texorpdfstring{$E_2$}{E2}}

The ring $K_*^\wedge E$ is a $K_*^\wedge K$-comodule. By \Cref{E cooperations}
(and the fact that $K_*^\wedge E$ is classically complete), it is equivalently a
continuous $\Z_p^\times$-module, and there's a homotopy fixed points spectral
sequence
\[ E_2 = H^*_{cts}(\Z_p^\times, K_*^\wedge E) \Rightarrow \pi_*\LE. \]

\begin{proposition}
  The cohomology of $\Z_p^\times$ acting on $K_*^\wedge E$ is concentrated in
  degree zero, so that the above spectral sequence collapses at $E_2$.
\end{proposition}

For the next proposition, write $\sigma$ for the Frobenius map of any ring or
scheme of characteristic $p$. If $\G$ is a formal group over a ring $R$ of
characteristic $p$, we define the formal group $\G^{(p)}$ and the
\textbf{relative Frobenius} $\Frob:\G \to \G^{(p)}$ by the following pullback
square:
\[ \xymatrix{ \G \ar[dr]^{\Frob} \ar@/^/[drr]^\sigma \ar@/_/[ddr] & & \\
    & \G^{(p)} \ar[r] \ar[d] \pullbackcorner & \G \ar[d] \\
    & \Spec R \ar[r]_\sigma & \Spec R } \]

\begin{proposition}
  Every $\theta$-algebra structure on $\LE_0$ extends to a $\theta$-comodule
  algebra structure on $K_0^\wedge E$.
\end{proposition}

\begin{proof}
  Fix a $\theta$-algebra structure on $\LE_0$. Since $K_0^\wedge E$ is
  torsion-free, it suffices to extend $\psi^p$ over $K_0^\wedge E$. Note that
  $\psi^p$ should not be $K_0$-linear in the most obvious sense: rather, $K_0 =
  Wk$ carries a unique Frobenius lift (which we will also write $\sigma$), and
  $\psi^p:K_0^\wedge E \to K_0^\wedge E$ should be $K_0$-linear where $K_0$ acts
  on the target via $\sigma$.

  By \Cref{FE theorem}, to define $\psi^p:K_0^\wedge E \to K_0^\wedge E$, it
  suffices to define a $j:\LE_0 \to K_0^\wedge E$ and an isomorphism of formal
  groups,
  \[
    \gamma: ( \Gamma_1 \otimes_k^i K_0^\wedge E/p )^{(p)} = \Gamma_1
    \otimes_{k}^{i\sigma} K_0^\wedge E/p \stackrel{\sim}{\to} \H
    \otimes_{k((x))}^{\ol{j}} K_0^\wedge E/p.
  \]
  Note that the twisted $K_0$-algebra structure on $K_0^\wedge E/p$ also forces
  a Frobenius twist on the map $i$ used to identify $\Gamma_1$ with the formal
  group of $K_0^\wedge E/p$.

  Define
  \[ j: \LE_0 \stackrel{\psi^p}{\to} \LE_0 \inj K_0^\wedge E. \] Observe that
  the reduction of $j$ mod $p$ is
  \[ \ol{j}:k((x)) \stackrel{\sigma}{\to} k((x)) \inj K_0^\wedge E/p, \] so that
  \[ \H \otimes_{k((x))}^{\ol{j}} K_0^\wedge E/p = (\H
    \otimes_{k((x))}^{\ol{\mathrm{can}}} K_0^\wedge E/p)^{(p)}, \] where we have
  written $\mathrm{can}$ for the canonical morphism $\LE_0 \inj K_0^\wedge E$.

  Now $K_0^\wedge E$ carries an isomorphism
  \[ \gamma_{ \mathrm{can} }: \Gamma \otimes_k K_0^\wedge E/p
    \stackrel{\sim}{\to} \H \otimes_{k((x))}^{\ol{j}_{\mathrm{can}}} K_0^\wedge
    E/p; \] we define $\gamma$ to be the base change of $\mathrm{can}$ along the
  Frobenius,
  \[ \gamma = \gamma_{ \mathrm{can} }^{(p)}: \Gamma_1 \otimes_k^{i\sigma}
    K_0^\wedge E/p \stackrel{\sim}{\to} \H
    \otimes_{k((x))}^{\ol{j}_{\mathrm{can}}\sigma} K_0^\wedge E/p.
  \]

  The given $(j, \gamma)$ define a map $\psi^p:K_0^\wedge E \to K_0^\wedge E$,
  which we have to show lifts the Frobenius. In other words, we must show that
  the square
  \[
    \xymatrix{ K_0^\wedge E \ar[r]^{\psi^p} \ar[d] & K_0^\wedge E \ar[d] \\
      K_0^\wedge E/p \ar[r]_{\sigma} & K_0^\wedge E/p }
  \]
  commutes. Both compositions are $K_0$-algebra maps, when $K_0$ acts on the
  target $K_0^\wedge E/p$ via the Frobenius lift $\sigma$, so we can check that
  they are equal using \Cref{FE theorem}. The top right composition represents
  the pair
  \[ (\ol{j} = \ol{\psi^p} = \sigma_{k((x))}, \, \gamma =
    \gamma_{\mathrm{can}}^{(p)}). \] Applying a Frobenius to $K_0^\wedge E/p$
  clearly twists its subring $k((x))$ by a Frobenius. Finally, $K_0^\wedge E$
  carries the coefficients of a power series which is the universal isomorphism
  $\gamma_{\mathrm{can}}$ between $\Gamma_1$ and $\H$; applying the Frobenius to
  the coefficients of these power series is the same as pulling back this
  isomorphism along the Frobenius.

  Finally, to show that we have defined a $\theta$-comodule algebra structure on
  $K_0^\wedge E$, we have to show that the $\psi^p$ just defined commutes with
  the action of the Morava stabilizer group. In general, a group element $(\tau,
  g) \in \Aut(k, \Gamma_1)$ acts on a point $(j, \gamma)\in
  \Hom_{cts}(F_0^\wedge E,R)$ by precomposition with $\gamma$:
  \[
    (\tau, g)(j, \gamma) = (j, \gamma g^{-1}:\Gamma_1 \otimes_{k}^{i\tau}
    R/\mathfrak{m} \stackrel{g^{-1}}{\to} \Gamma_1 \otimes_k^i R/\mathfrak{m}
    \stackrel{\gamma}{\to} R/\mathfrak{m}).
  \]
  This operation commutes with pullback of $\gamma$ along the Frobenius.
\end{proof}

\subsection{Obstruction theory}

In \cite{GH1, GH2}, Goerss and Hopkins construct an obstruction theory for
realizing a $\theta$-comodule algebra as the completed $K$-homology of a
$K(1)$-local $\Einf$ ring spectrum. Before describing the obstruction theory, we
must make some preliminary definitions. Note that our situation is slightly
simplified by the fact that $K_0^\wedge E$ is classically $p$-complete (rather
than just L-complete) and concentrated in even degrees; the definitions below
need to be modified to handle more general cases.

\begin{definition}[{\cite[Definition 2.2.7]{GH2}}] Let $A$ be a $p$-complete
  $\theta$-comodule algebra. An \textbf{$A$-$\theta$-module} is a $p$-complete
  $A$-module $M$ equipped with a continuous action by $\psi^k$, $k \in
  \Z_p^\times$, and an operation $\theta:M \to M$, such that
  \begin{itemize}
  \item for $a \in A$ and $m \in M$, $\psi^k(am) = \psi^k(a)\psi^k(m)$;
  \item for $a \in A$ and $m \in M$, $\theta(am) = a^p\theta(m) +
    p\theta(a)\theta(m)$.
  \end{itemize}
  Write $\Mod_A^\theta$ for the category of $A$-$\theta$-modules.
\end{definition}

\begin{remark}\label{split square-zero}
  There is an equivalence between $\Mod_A^\theta$ and the category of abelian
  group objects in $( \ComodAlg_\theta )_{/A}$, given by sending
  \[ [B \to A] \in \Ab((\ComodAlg_\theta)_{/A}) \,\, \mapsto \,\, \ker(B \to
    A).
  \]
  Indeed, as a ring, $B$ is a split square-zero extension of $A$, of the form
  $A \oplus M \stackrel{\mathrm{proj}}{\to} A$; the structure on $M$ induced
  by the $\theta$-comodule algebra structure on $A$ is precisely the
  $A$-$\theta$-module structure defined above.
\end{remark}

\begin{definition}
  Let $A$ be a $\theta$-comodule algebra, and $M$ an $A$-$\theta$-module. A
  \textbf{derivation} from $A$ into $M$ is a map $A \to A \oplus M$ of
  $\theta$-algebras augmented over $A$, where $A \oplus M$ is the split
  square-zero extension defined in \Cref{split square-zero}. We write
  $\Der_{\ComodAlg_\theta}(A, M)$ for the abelian group of derivations from
  $A$ into $M$.
\end{definition}

\begin{definition}
  The \textbf{Andr\'e-Quillen cohomology groups} of the $\theta$-algebra $A$,
  with coefficients in $A$-$\theta$-modules $M$, are the left derived functors
  \[ D^s_{\ComodAlg_\theta}(A, \cdot) \] of $\Der_{\ComodAlg_\theta}(A, M)$.
\end{definition}

\begin{remark}
  Making sense of these derived functors requires defining a well-behaved
  model structure on simplicial $\theta$-algebras, which is done in
  \cite[Proposition 2.3.1]{GH2}.
\end{remark}

\begin{remark}
  Completely analogously but forgetting about the $\Z_p^\times$-action
  everywhere, we can define $\theta$-modules over a $\theta$-algebras,
  derivations, and Andr\'e-Quillen cohomology.
\end{remark}

We can now state the obstruction theory. This is a special case of
\cite[Theorem 3.3.7]{GH2} as well as the results on $\theta$-algebras in
\cite{GH2}; see \cite{Szymik} for a more precise statement.

\begin{thm}\label{height 1 obstruction theory}
  Let $A$ be an $p$-complete $\theta$-comodule algebra. Then there are
  successively defined obstructions to realizing $A$ as $K_0^\wedge X$, where
  $X$ is an $\Einf$ algebra such that $K_*^\wedge X$ is concentrated in even
  degrees, in the Andr\'e-Quillen cohomology groups
  \[ D^{s+2}_{\Comod\Alg_\theta}(A, \Omega^{s} A), \,\, s \ge 1. \] There are
  successively defined obstructions to the uniqueness of this realization in
  \[ D^{s+1}_{\Comod\Alg_\theta}(A, \Omega^{s} A), \,\, s \ge 1. \]
\end{thm}

We now prove that any $\theta$-algebra structure on $\CLaur$ is induced by an
$\Einf$-algebra structure on $\LE$, using the argument of \cite[Section
2.4.3]{GH2}. Thanks again to Dominik Absmeier for calling my attention to a
mistake in the earlier version of this statement.

\begin{thm}\label{LK(1)E2}
  For any $p$-complete $\theta$-comodule algebra $A$ such that $A$ is
  isomorphic to $K_0^\wedge E$ as a ring with $\Z_p^\times$-action, there is
  an even periodic $\Einf$-algebra $X$ with $K_0^\wedge X = A$ as
  $\theta$-comodule algebras, and with $X \simeq L_{K(1)}E_2$ as homotopy commutative
  ring spectra.
\end{thm}

\begin{proof}
  We want to show that the obstruction groups
  \[ D^{s+2}_{\Comod\Alg_\theta}(A, \Omega^s A) \] vanish for $s \ge 1$. Since
  $\Omega^s A \cong K_0^\wedge \Omega^s E$ is an extended $\Z_p^\times$-module
  by \Cref{action extended}, these reduce \cite[Proposition 2.4.7]{GH2} to
  Andr\'e-Quillen cohomology of $\theta$-algebras without
  $\Z_p^\times$-action:
  \[ D^{s+2}_{\Comod\Alg_\theta}(A, \Omega^s A) \cong D^{s+2}_{\Alg_\theta}(A,
    \Omega^s \LE_0). \] The complete cotangent complex of $A$ over $Wk$ is a
  $\theta$-module, and there is a composite functor spectral sequence
  \[ \Ext^p_{\Mod_{\theta,A/}}(\pi_q\L_{A/Wk}, \Omega^s \LE_0) \Rightarrow
    D^{p+q}(A, \Omega^s \LE_0). \] But $Wk = K_0 \to K_0^\wedge E = A$ is
  formally smooth by \Cref{FE theorem} and \Cref{formally smooth}. Thus,
  $\L_{A/Wk}$ is just the K\"ahler differentials $\Omega_{A/Wk}$ concentrated
  in degree zero, and these are a $p$-completion of a free module. Finally,
  there is a resolution $\Omega_{A/Wk}$ by free $\theta$-modules over $A$,
  \[ 0 \to A[\theta] \otimes_A \Omega_{A/Wk} \stackrel{\theta}{\to} A[\theta]
    \otimes_A \Omega_{A/Wk} \to \Omega_{A/Wk}. \] For any complete
  $\theta$-module $M$ over $A$,
  \[ \Ext^*_{\Mod_{\theta,A/}}(A[\theta] \otimes_A \Omega_{A/Wk}, M) =
    \Ext^*_{\Mod_A}(\Omega_{A/Wk}, M), \] which is concentrated in degree zero
  because $A$ is pro-free. Thus, the Andr\'e-Quillen cohomology groups are
  concentrated in cohomological degrees 0 and 1, and in particular, those that
  can contain obstructions vanish.

  This produces a $K(1)$-local $\Einf$-algebra $X$ with $K_0^\wedge X = A$.
  Since $K_0^\wedge X \cong K_0^\wedge E$ as $\Z_p^\times$-modules, we also have
  $K_t^\wedge X \cong K_t^\wedge E$ as $\Z_p^\times$-modules for all $t$. In
  particular, these are extended and $\Z_p$-free, and so the Adams-Novikov
  spectral sequence
  \[
    E_2 = \Ext_{\Z_p^\times}^*(K_*^\wedge E, K_*^\wedge X) \Rightarrow [E,
    X]_*
  \]
  is concentrated on the 0-line. Thus the isomorphism $K_*^\wedge E \to
  K_*^\wedge X$ of comodule algebras lifts to an equivalence $E \to X$ of
  homotopy commutative ring spectra. 
\end{proof}

\begin{corollary}\label{LK(1)E2 nonisomorphic}
  There exist $\Einf$-algebras with underlying ring spectrum $L_{K(1)}E_2$,
  which are not equivalent to the $K(1)$-localization of the $\Einf$-algebra
  $E_2$.
\end{corollary}

\begin{proof}
  This follows from \Cref{LK(1)E2} and \Cref{nonisomorphic theta}.
\end{proof}

\end{document}